\input ./style/arxiv-general.cfg
\documentclass[aos,MSNbibl,nameyear,secthm,seceqn,dvips]{arximspdf}
\makeatletter
   \@ifpackageloaded{graphicx}{}{\usepackage{graphicx}}
\makeatother
\usepackage{algorithm,algorithmic}

\doi{10.1214/14-AOS1290}
\volume{43}
\issue{3}
\pubyear{2015}
\firstpage{1027}
\lastpage{1059}
\docsubty{FLA}

\makeatletter
\newcommand{\ds}{\displaystyle}
\newcommand{\lleft}{\left}
\newcommand{\rright}{\right}
\newcommand{\rrVert}{\Vert}
\newcommand{\llVert}{\Vert}
\newtheorem{lemma}[thm]{Lemma}
\newtheorem{corollary}[thm]{Corollary}
\newproclaim{remark}{Remark}[section]
\makeatother

\begin{document}
\begin{frontmatter}

\title{Robust and computationally feasible community detection in the
presence of arbitrary outlier~nodes\thanksref{T1}}
\runtitle{Robust community detection}

\begin{aug}
\author[A]{\fnms{T. Tony}~\snm{Cai}\ead
[label=e1]{tcai@wharton.upenn.edu}}
\and
\author[A]{\fnms{Xiaodong}~\snm{Li}\ead[label=e2]{xiaodli@wharton.upenn.edu}\corref{}}
\runauthor{T. T. Cai and X. Li}
\affiliation{University of Pennsylvania}
\address[A]{Department of Statistics \\
The Wharton School\\
University of Pennsylvania\\
400 Jon M. Huntsman Hall\\
3730 Walnut Street\\
Philadelphia, Pennsylvania 19104-6340\\
USA\\
\printead{e1}\\
\phantom{E-mail:\ }\printead*{e2}}
\end{aug}
\thankstext{T1}{Supported in part by NSF Grants
DMS-12-08982 and DMS-14-03708, NIH Grant R01 CA127334-05, and the Wharton
Dean's Fund for Post-Doctoral Research.}

%
\received{\smonth{4} \syear{2014}}
%
\revised{\smonth{11} \syear{2014}}

%
\begin{abstract}
Community detection, which aims to cluster $N$ nodes in a given graph
into $r$ distinct groups based on the observed undirected edges, is an
important problem in network data analysis. In this paper, the popular
stochastic block model (SBM) is extended to the generalized stochastic
block model (GSBM) that allows for adversarial outlier nodes, which are
connected with the other nodes in the graph in an arbitrary way. Under
this model, we introduce a procedure using convex optimization followed
by $k$-means algorithm with $k=r$.

Both theoretical and numerical properties of the method are analyzed.
A~theoretical guarantee is given for the procedure to accurately detect
the communities with small misclassification rate under the setting
where the~number of clusters can grow with $N$. This theoretical result
admits to the best-known result in the literature of computationally
feasible community detection in SBM without outliers. Numerical results
show that our method is both computationally fast and robust to
different kinds of outliers, while some popular computationally fast
community detection algorithms, such as spectral clustering applied to
adjacency matrices or graph Laplacians, may fail to retrieve the major
clusters due to a small portion of outliers. We apply a slight
modification of our method to a political blogs data set, showing that
our method is competent in practice and comparable to existing
computationally feasible methods in the literature. To the best of the
authors' knowledge, our result is the first in the literature in terms
of clustering communities with fast growing numbers under the GSBM
where a portion of arbitrary outlier nodes exist.
\end{abstract}

%
\end{frontmatter}

\section{Introduction}
\label{intro.sec}

Driven by applications in a wide range of fields, including
engineering, genomics, sociology, psychology and computer science,
analysis of graph and network data has drawn significant recent
interest. Random graph models have been introduced to characterize the
structure of the networks and a large number of algorithmic approaches
have been proposed for various applications. See, for example,
\citeauthor{Fienberg2010} (\citeyear{Fienberg2010,Fienberg2012}), \citet{GZFA2010}, and the references therein for
overviews and recent work.

An important problem in the analysis of network data is that of
community detection which aims to cluster the nodes in a given graph
into distinct groups or communities based on the observed undirected
edges. Community detection has proven to be both technically and
computationally challenging. It also has deep connections to other
fields such as spin-glass theory and signal processing. In terms of
statistical modeling, the most well-known model for community detection
is perhaps the stochastic block model (SBM) proposed in \citet{HLL1983}.
Under the SBM, the graph of interest is assumed to be a random one with
independent edges, and the within-group edge density is assumed to be
greater than the between-group edge density.

To be specific, suppose $G=(V, E)$ is a random graph where $V$ is a
fixed set of vertices consisting of $n$ nodes, and $E$ is a random set
of edges. Assume that the $n$ nodes are indexed by $[n]:=\{1, \ldots,
n\}$ and each of these nodes belongs to one and only one of the $r$
nonoverlapping groups. This amounts to assigning each node $j \in[n]$
a group label by a labeling function $\phi(j)\in\{1, \ldots, r\}$. We
denote by $\mathbf{A}=(A_{ij})_{1\leq i, j \leq n}$ the random adjacency
matrix of this random graph. Then for each pair $(i,j)$, $1\leq i, j
\leq n$, $A_{ij}=0\mbox{ or }1$, indicating whether the nodes $i$ and
$j$ are connected or not, respectively. We only consider undirected
graph with no self loops, so $\mathbf{A}$ is symmetric, and all its
diagonal entries are $0$. For pairs $(i,j)$ with $1\leq i< j \leq n$,
$A_{ij}$'s are assumed to be independent Bernoulli random variables
with parameters $B_{\phi(i)\phi(j)}$, where the symmetric matrix
$\mathbf{B} \in\mathbb{R}^{r\times r}$ is referred to as the \textit{connectivity} matrix. In a basic model, denote by $q^+$ and $p^-$ the
maximum cross-group density and the minimum within-group density, namely
%
\begin{equation}
\label{eqdefp}
q^+:=\max_{1\leq i<j\leq r} B_{i j}, \qquad  p^-:=\min
_{1\leq i\leq r} B_{ii}.
\end{equation}
Moreover, the within-group densities are assumed to be greater than the
cross-group densities, that is,
%
\begin{equation}
\label{eqdefdelta}
p^- - q^+:=\delta>0.
\end{equation}
This is a common assumption in the literature of community detection
under the SBM; see, for example, \citet{RCY2011,CCT2012}. Denote the
minimum community size by $n_{\min}:= \min_{1\leq l\leq r}
|\phi
^{-1}(l)|$, where $|S|$ denotes the cardinality of the set $S$. Then
the difficulty of the community detection problem is determined by the
tuple $(n, r, q^+,  p^-, n_{\min})$.


%

Under the SBM, various community detection algorithms have been
proposed and studied in the literature, with different emphases on
computational complexity and statistical accuracy. These include greedy
algorithms, such as hierarchical agglomeration [see, e.g., \citet
{CNM2004}]; greedy \mbox{methods} guided by global criterion maximization,
such as modularity function maximization [see, e.g., \citet{NG2004}] and
profile likelihood function maximization [see, e.g., \citet{BC2009,ZLZ2012}]; stochastic model based methods, such as variational
likelihood methods [see, e.g., \citet{BCCZ2013,CDP2012}],
pseudo-likelihood methods with EM algorithm [see, e.g., \citet
{ACBL2013}], Bayesian methods with Gibbs sampling, Markov chain Monte
Carlo and belief propagation [see, e.g., \citet{SN1997,NS2001,DKMZ2011}]; graph distance methods [see, e.g., \citet{BB2013}]; spectral
clustering, its variations and other spectral methods [see, e.g., \citet
{McSherry2001,GM2005,RCY2011,CCT2012,CL2009,BXKS2011,STFP2012,FSTVP2013,Jin2012,JY2013,SB2013,LR2013}]; and convex optimization
methods [see, e.g., \citet{MS2010,OH2011,JCSX2011,AV2011,CSX2013,Ames2013}].

Among these methods, greedy methods are usually computationally
feasible, while their statistical accuracy has not been fully
established in theory. Modularity or profile likelihood methods are
proven to be consistent when the number of groups is fixed. However,
they are in principle computationally NP hard. Similarly, stochastic
model based methods are usually computationally difficult and not fully
justified in theory. Spectral clustering is a popular algorithm for
community detection, since it is fast in computation and easy to
implement. It has been proven that spectral clustering is consistent
even when the number of groups $r$ grows on the order of $O(\sqrt{n})$.
Although in practice spectral clustering is believed to work well only
for dense graphs, several recent papers, \citet{ACBL2013,SB2013,JY2013,LR2013}, have shown that spectral clustering or its variations
also work well for sparse graphs.

The SBM is admittedly an oversimplified model for many applications,
and different generalizations have been proposed in the literature,
which encompass mixture model [see \citet{NL2007}], where the parametric
model for the connectivity probabilities is based on the relationship
between vertices and groups, instead of between different groups;
degree corrected model [see \citet{CL2009,KN2011,ZLZ2012}]; Latent
variable method [see \citet{HRT2007}] and mixed membership model [see
\citet{ABFX2008}]. However, each of these GSBMs focuses on a single
latent graph structure, while in practice, due to lack of information,
this additional structure is not easy to detect if it only applies to a
few nodes of the graph. Different types of outliers may appear in a
single graph, and it is difficult to use a complex generalization of
the SBM to model multiple types of outlier nodes. The SBM is usually
the first model to fit the data because of its simple form, even if it
is believed that there is possibly a small portion of nodes which are
not modeled well.
Robustness in presence of arbitrary outliers is an important property
for given community detection algorithms. In this paper we consider
robust community detection in the presence of arbitrary outlier nodes,
and the main question we wish to answer is the following:

\begin{quote}
Does there exist a computationally fast community detection
method that is robust to a portion of arbitrary outlier nodes with
theoretical guarantees?
\end{quote}

Our answer is affirmative, and we will introduce our model,
methodology, numerical results and theoretical guarantees with rigorous
proofs in this paper. We begin by formalizing the GSBM which allows for
a small portion of arbitrary nodes.

\subsection{Generalized stochastic block model}

We introduce a flexible model for community detection which covers a
range of settings in practice where the usual SBM is not suitable.
More specifically, we assume the undirected graph $G=(V, E)$ has
$N:=n+m$ nodes, among which there are $n$ ``inliers'' obeying the SBM
described above, while the other $m$ nodes are ``outliers'' which are
connected with the other nodes in an arbitrary way. We refer to this
model as \textit{generalized stochastic block model} (\textit{GSBM}). Denote
$V=[N]=I \cup O$, where $I$ is the set of indices of the inliers, while
$O$ is the set of indices of outliers. Each inlier node $i \in I$ is
assigned a label $\phi(i)\in\{1,\ldots, r\}$, while all outliers are
simply labeled $\phi(i)=r+1$. For any two nodes $i, j\in I$, $\mathbb{P}
((i,j)\in E)=B_{\phi(i)\phi(j)}$, and moreover we assume the event
$\{(i,j)\in E\}$, $i < j \in I$ are independent. The $r\times r$ symmetric
connectivity matrix $\mathbf{B}$ only represents the likelihood of
connectivity of the inlier nodes. The connectivity between the outliers
and the inliers and the connectivity among the outliers themselves are
arbitrary. The only restriction of the connectivity of the outliers is
that there is no self-loop.

The GSBM can be equivalently expressed in terms of its adjacency matrix
$\mathbf{A}$. To be specific, define
%
\begin{equation}
\label{eqadjacency}
\mathbf{A}=\mathbf{P}
\lleft[\matrix{ \mathbf{K} &
\mathbf{Z} \vspace*{3pt}
\cr
\mathbf {Z}^{\intercal} & \mathbf{W} } \rright]
\mathbf{P}^{\intercal}=\mathbf{P}
\lleft[
\matrix{ \mathbf{K}_{11} &\ldots& \mathbf{K}_{1r} & \mathbf{Z}_1\vspace*{3pt}
\cr
\vdots&\ddots& \vdots& \vdots\vspace*{3pt}
\cr
\mathbf {K}_{1r}^{\intercal} & \ldots& \mathbf{K}_{rr} &
\mathbf{Z}_r \vspace*{3pt}
\cr
\mathbf {Z}_1^{\intercal}
& \cdots& \mathbf{Z}_r^{\intercal} & \mathbf{W} } \rright]
\mathbf{P}^{\intercal},
\end{equation}
where $\mathbf{W}\in\mathbb{R}^{m\times m}$ is an arbitrary symmetric
$0$--$1$ matrix with all diagonal entries being $0$, $\mathbf{Z} \in
\mathbb
{R}^{n\times m}$ is an arbitrary $0$--$1$ matrix, $\mathbf{P}$ is an
unknown $N\times N$ permutation matrix, in which there is only one $1$
in each row and column, while all other entries are $0$'s, and $\mathbf
{K}$ is an $n\times n$ symmetric matrix which captures the connectivity
of the inliers, thus corresponding to the usual SBM. The off-diagonal
entries of $\mathbf{K}$ are independent Bernoulli variables, with
parameter $B_{ij}$ if the entry belongs to the submatrix $\mathbf
{K}_{ij}$. Denote the dimension of $\mathbf{K}_{ii}$ to be $l_i$ for $i=1,
\ldots, r$. Then $n=\sum_{i=1}^r l_i$. Similar to SBM, $n_{\min
}=\min_{1\leq i\leq r}l_i$. The parameters $p^-$ and $q^+$ are defined as in
(\ref{eqdefp}) and $\delta$ in (\ref{eqdefdelta}).
Then the difficulty of community detection under the GSBM is
parameterized by the tuple $(n, m, r, p^-, q^+, n_{\min})$.

Here we emphasize that $\mathbf{Z}$ and $\mathbf{W}$ are not necessarily
fixed with respect to the randomness of $\mathbf{K}$. Both $\mathbf
{Z}$ and
$\mathbf{W}$ can depend on $\mathbf{K}$ in arbitrary forms. In other words,
the connectivity between the outliers and the inliers is allowed to
depend on the connectivity among the inlier nodes. This is also a
generalization of standard SBM, where the connectivity between each
pair of nodes is stochastically independent of the connectivity between
other pairs.

The GSBM is a flexible model and is widely applicable. It covers
various types of outliers which are common in practice, and we name a
few as follows:
\begin{itemize}
\item \textit{Mixed membership}. The SBM assumes that each node belongs to
one and only one predetermined cluster. If most nodes obey this
property, while there is a small portion of nodes each belonging to
more than one clusters, these nodes are referred to as having mixed membership.
When only a small portion of nodes have mixed membership, it is natural
to treat them as outliers in an ordinary SBM.

\item \textit{Hubs}. In social networks and others, it is natural that
some nodes have many more connections than most of others. Moreover, it
is possible that these nodes belong to several groups without obvious
bias to any specific one. These nodes are referred to as hubs, and can
be treated as outliers in our GSBM.

\item \textit{Small clusters}. The SBMs are usually employed to model big
and significant clusters, while small clusters are difficult to detect.
Small clusters are often not detectable because they are too small and
possibly weak.
The number of small clusters is also difficult to estimate; however,
this information is essential for most popular algorithms in the
literature, such as spectral clustering and modularity methods. 
The nodes in the small clusters can be treated as outliers in our GSBM.

\item \textit{Independent neutral nodes}.  In a given graph, in addition to
the well-classified nodes, there might be some nodes which do not
belong to any significant groups, and also have fewer connections than
most other nodes. We refer to these objects as independent neutral
nodes. For example, in the political blogs data set introduced later, a
small portion of blogs have very few connections. Such blogs may have
strong preference in politics; however, this cannot be seen from only
the graph representation. Therefore, these nodes are regarded as
independent neutral nodes, which are naturally taken as outliers.
\end{itemize}

In practice, it is difficult or even impossible to modify the usual SBM
to model precisely the possible combinations of mixed membership, hubs,
small clusters, independent neutral nodes and other types of settings.
Moreover, complex statistical models may also result in overfitting and
high computational complexity in clustering. Therefore, the SBM is
usually set up based on the basic properties of the graph. For example,
in the political blogs network application discussed in Section~\ref{application.sec}, an SBM with $2$ clusters is preferred, since it is
known that there are mainly two significant clusters: liberals and
conservatives. However, it is also known that there are many
independent groups advocating various causes that lie outside of the
two main clusters.

The GSBM can also be taken as a criterion to evaluate the robustness of
community detection algorithms. When an SBM is adopted based on the
properties of most nodes of a given graph, or equivalently, most nodes
can be well modeled by an SBM in use, the robustness of a given
community detection algorithm depends on whether a small portion of
outliers will completely change the clustering result, or most nodes
can be still well clustered. Therefore, a graph clustering algorithm is
robust if it is guaranteed to have good performance under the GSBM.

\subsection{Organization of the paper}
The rest of the paper is organized as follows. In Section~\ref{method.sec} the method of convexified likelihood method is introduced,
followed by a detailed alternating directional augmented Lagrangian
algorithm. Section~\ref{theory.sec} is focused on the theoretical
consistency of the convex optimization method in the inference of the
underlying groups specified by the GSBM. Numerical results on the
analyses of the simulated data and a real data set about political
blogs are presented in Section~\ref{simulation.sec}. A discussion is
given in Section~\ref{discussion.sec}, and the proofs of the main
theoretical results are contained in Section~\ref{proof.sec}.
Additional technical proofs are given in the Supplementary Material.

\section{Methodology}
\label{method.sec}

In this section we propose a community detection algorithm which is
robust and computationally feasible with theoretical guarantee of
consistency. In the literature, greedy algorithms such as hierarchical
clustering are not fully justified in theory, while modularity and
profile maximum likelihood methods are computationally NP hard.
Stochastic model based methods, such as maximum likelihood or
variational likelihood method, have been proven to have certain
consistency when the number of blocks is fixed as the number of nodes
going to infinity. However, they are also computationally difficult. EM
algorithm is naturally proposed for solving relevant maximum likelihood
formulation, but there is no theoretical guarantee of convergence with
reasonable rate. Bayesian methods such as Gibbs sampling and belief
propagation have also been proposed in the literature without rigorous
theoretical justifications.

Unlike the aforementioned methods, the spectral clustering methods have
the advantage of fast algorithms. Spectral clustering algorithms are
easy to implement because there is no tuning parameter. Moreover,
strong theoretical results have been established under various
conditions; see the references mentioned in the previous section.
However, as indicated in \citet{JY2013}, ordinary spectral clustering
applied to the graph Laplacian may not work due to the existence of
small and weak clusters. We use a simulated data set to illustrate that
ordinary spectral clustering applied to the graph Laplacian or the
adjacency matrix is not consistent under the GSBM. Other types of
numerical examples can be found in \citet{JY2013}.

\begin{figure}

\includegraphics{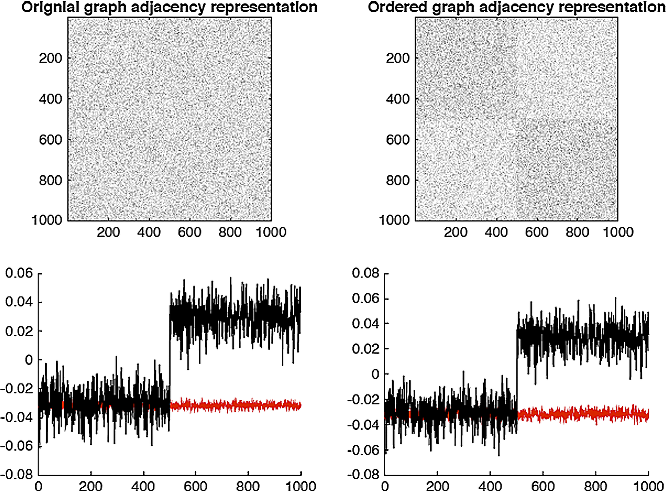}

\caption{The upper left panel illustrates the adjacency matrix of
$1000$ nodes satisfying the ordinary SBM. The upper right panel is the
adjacency matrix obtained by permuting the adjacency matrix such that
nodes $1$ to $500$ belong to the same cluster while the remaining ones
constitute another cluster. The lower left panel plots the eigenvectors
of the graph Laplacian corresponding to the top $2$ eigenvalues in
absolute value (red for the first and black for the second), while
those for the adjacency matrix are plotted in the lower right panel. In
both cases, these two eigenvectors are capable of discriminating
between the two communities.}
\label{figbig2}
\end{figure}

First, we create a data set of $n=1000$ nodes obeying the ordinary SBM
with $r=2$ clusters. We also assume that the two clusters are perfectly
balanced; that is, there are $500$ nodes in each cluster. The
within-group probability is $p=0.17$, while the cross-group probability
is $q=0.11$. Under this set-up, the adjacency matrix is shown as in
Figure~\ref{figbig2}.

\begin{figure}

\includegraphics{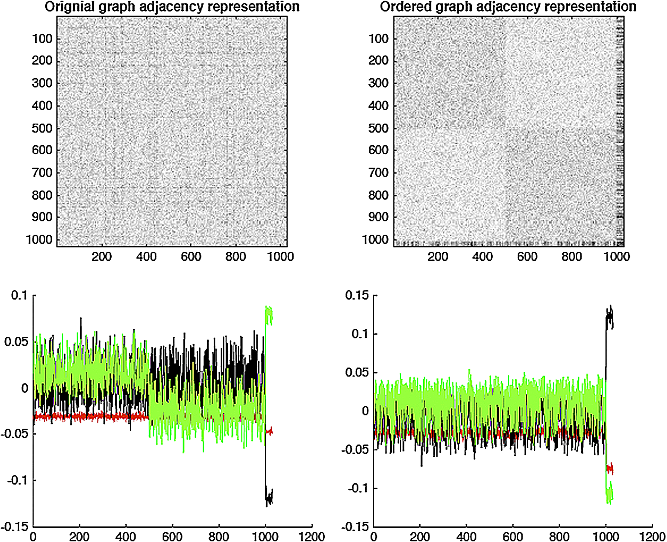}

\caption{The upper left panel illustrates the adjacency matrix of
$1030$ nodes satisfying the GSBM with two major clusters and $30$
outliers. The upper right panel is obtained by permuting the nodes such
that nodes belonging to the same group are consecutive. The lower left
panel plots the eigenvectors of the graph Laplacian corresponding the
top $3$ eigenvalues in absolute value (red for the first, black for the
second and green for the third), while those for the adjacency matrix
are plotted in the lower right panel. Ordinary spectral clustering with
$r=2$ or $r=3$ is ineffective or even powerless on this data set since
the top three eigenvectors cannot clearly discriminate between the two
main communities.}
\label{figbig2small2}
\end{figure}

Spectral clustering applied directly to either the graph Laplacian or
the adjacency matrix of this graph data has good performance of
clustering. To illustrate this, we plot the eigenvectors corresponding
to the top two eigenvalues (in absolute value) of the graph Laplacian
and the adjacency matrix, respectively, in Figure~\ref{figbig2}. In
each case, the two eigenvectors combined are capable of discriminating
between the two clusters. Therefore, spectral clustering methods work
for our data set when there are no outliers.

Now we consider the GSBM by adding only $m=30$ outliers into the above
model with $r=2$ clusters. To specify this GSBM, it suffices to explain
$\mathbf{Z}$~and~$\mathbf{W}$ in~(\ref{eqadjacency}). We assume
$\mathbf{W}$ is
the adjacency matrix of a random graph with $30$ nodes and independent
edges. Moreover, we assume the probability of connectivity is~$0.7$. We
define $\mathbf{Z}$ as a $1000 \times30$ with independent Bernoulli
entries. We also let $\mathbb{E}\mathbf{Z} = \bolds{\beta}\mathbf
{1}^{\intercal}=[\bolds{\beta}, \bolds{\beta}, \ldots, \bolds
{\beta}]$. The components of $\bolds{\beta}$ are $1000$ i.i.d.
copies of $U^2$, and $U$ is a uniform random
variable on $[0, 1]$. The unordered and ordered adjacency matrices are
given in Figure~\ref{figbig2small2}.

Suppose the data set is still modeled approximately by the SBM with
$r=2$. For this new data set, the two eigenvectors corresponding to the
top two eigenvalues (in absolute value) of the graph Laplacian or the
adjacency matrix cannot discriminate between the two major clusters.
Even if we treat the $30$ outliers as a single group due to their
homogeneous behavior in the graph and thereby use $r=3$, the third
eigenvector of the adjacency matrix is still unable to distinguish the
two major clusters. The third eigenvector of the graph Laplacian can
only discriminate a part of nodes in the two major clusters. Actually,
our numerical simulation shows that after applying spectral clustering
on the graph Laplacian with $r=3$, the misclassification rate among the
inliers is above $30$ percent. These three eigenvectors are plotted in
Figure~\ref{figbig2small2} for both cases. The figures indicate
that standard spectral clustering is not a robust community detection
method in the presence of very few adversarial outliers.

It was shown in \citet{JY2013} that under certain conditions, penalized
spectral clustering may reduce the effects of the small weak clusters,
but it is not clear whether penalized spectral clustering applied to
the graph Laplacian can diminish the influence of other types of
outliers. Another method to improve standard spectral clustering
methods is to detect outlier nodes based on the first several
eigenvectors. However, it is not clear whether there exists an approach
which can uniformly detect all kinds of outliers with a theoretical guarantee.

In order to find in one shot the major clusters among the inlier nodes,
we introduce in Section~\ref{procedure.sec} a convex optimization
method as well as a detailed algorithm which is implementable. It will
be shown in Section~\ref{theory.sec} that the proposed procedure is
robust against a small portion of arbitrary outliers with theoretical
guarantees.

\subsection{Convex optimization}
\label{procedure.sec}

In this section, we will choose the method of semidefinite programming
(SDP) to fit the GSBM, followed by a $k$-means clustering. Numerically,
SDP is well known to be computationally feasible, and various efficient
algorithms were proposed for solving different types of SDP.
Theoretically, under the ordinary SBM, SDP methods are shown to be
capable in detecting communities; see \citet{MS2010,OH2011,JCSX2011,AV2011,CSX2013}. We propose a new convex optimization method inspired
by existing SDP methods in the literature. The significantly novel part
is that we will prove that this SDP method can consistently cluster the
nodes when there is a portion of arbitrary type of outliers. The formal
statement is given in Section~\ref{theory.sec}, and all the proofs are
deferred to Section~\ref{proof.sec} and the supplemental article \citet
{CL2014}.

First, we derive the convex optimization from the viewpoint of fitting
a parametric model. This viewpoint was originally proposed in \citet
{CSX2013}, but we are going to derive a different convex optimization.
For now we only consider the ordinary SBM, which implies that $m=0$ and
$N=n$. By the definition of SBM, for all $1\leq i<j \leq n$, the events
$\{A_{ij}=1\}$ are independent. Recall that here $\mathbf{A}$ is the
observed adjacency matrix. Moreover, we define a symmetric matrix
$\mathbf{X}$ with all diagonal entries equal to $1$. For any $1\leq
i<j \leq
n$, we let $X_{ij}=0$ if the labeling functions $\phi(i)\neq\phi(j)$,
while $X_{ij}=1$ if $\phi(i)=\phi(j)$. Obviously, this matrix
$\mathbf{X}$
is of rank $r$ since there are $r$ groups.

Moreover, we consider a special case of the ordinary SBM. Suppose
$1>p>q>0$. For any $1\leq i<j \leq N$, when $X_{ij}=0$ let $\mathbb{P}
(A_{ij}=1)=q$; otherwise let $\mathbb{P}(A_{ij}=1)=p$. This gives
\[
\log\mathbb{P}(A_{ij}=1|X_{ij})=X_{ij}\log
p+(1-X_{ij})\log q
\]
and
\[
\log\mathbb{P}(A_{ij}=0|X_{ij})=X_{ij}
\log(1-p)+(1-X_{ij})\log(1-q).
\]
Since $\{A_{ij}=1\}$ are independent events, we have the log-likelihood function
\begin{eqnarray*}
\ell(\mathbf{A}|\mathbf{X}) &=& \sum_{1\leq i<j\leq n}  \bigl[A_{ij}\bigl(X_{ij}\log p+(1-X_{ij})\log q
\bigr)
\\
&&\hspace*{28pt}\quad{}+(1-A_{ij}) \bigl(X_{ij}\log(1-p)+(1-X_{ij})
\log(1-q)\bigr) \bigr].
\end{eqnarray*}
For any fixed $p$ and $q$, given $\mathbf{A}$, we would like to choose an
appropriate $\mathbf{X}$ to maximize $\ell(\mathbf{A}|\mathbf{X})$.
If we let
%
\begin{equation}
\label{eqlambdapq}
\lambda=\frac{\log(1-q)-\log p}{\log p-\log q+\log(1-q)-\log(1-p)},
\end{equation}
since the diagonal entries of $\mathbf{A}$ are all equal to $0$, the
maximization is equivalent to
\[
\max_{\mathbf{X}} \bigl\langle\mathbf{X}, (1-\lambda )\mathbf{A}-\lambda
(\mathbf{J}_N -\mathbf{I}_N-\mathbf{A}) \bigr\rangle,
\]
where $\mathbf{J}_N$ is the $N\times N$ matrix with all entries $1$. Now
let us figure out the constraint of $\mathbf{X}$. By the SBM, it is easy
to check that $\mathbf{X}$ must have the following form:
%
\begin{equation}
\label{eqmaximization}
\mathbf{X}=\mathbf{P}
\lleft[\matrix{
\mathbf{J}_{l_1} & & \vspace*{3pt}
\cr
&\ddots& \vspace*{3pt}
\cr
& &
\mathbf{J}_{l_r} } \rright] \mathbf{P}^{\intercal},
\end{equation}
where $\mathbf{P}$ is some unknown permutation matrix, while $\mathbf{J}_s$
is an $s\times s$ matrix with all entries $1$'s. Solving optimization
(\ref{eqmaximization}) under such constraint is computationally
infeasible, so we seek for some relaxed form. Here we notice there are
three major features of $\mathbf{X}$. First, it is positive semidefinite;
second, all its entries are between $0$ and $1$; third, it is of
rank-$r$, which is relatively low. If we convexify the second integer
constraint and neglect the third requirement, the relaxed maximum
likelihood method becomes
\begin{eqnarray*}
&& \max \qquad \bigl\langle\widetilde{\mathbf{X}}, (1-\lambda )\mathbf{A}-
\lambda(\mathbf{J}_N -\mathbf{I}_N-\mathbf{A}) \bigr
\rangle
\\
&&\mbox{subject to} \qquad \widetilde{\mathbf{X}}\succeq\mathbf{0},
\\
&& \hspace*{21pt}\qquad\qquad 0\leq\widetilde{X}_{ij}\leq 1 \qquad \mbox{for } 1\leq i, j\leq N.
\end{eqnarray*}
The above optimization method is different from that in \citet{CSX2013},
where the relaxation is based on the observation that $\mathbf{X}$ is of
low rank and hence a nuclear norm penalization is added up to the
original objective function. On the contrary, our convex relaxation is
derived from the observation that $\mathbf{X}$ is both low-rank and
positive semidefinite, and consequently we impose constraint of the
positive semidefinite cone.

Now let us come back to the robust community detection under the GSBM.
To control the possible outliers as formalized in the GSBM model, for
the convenience of theoretical analysis, we add an additional term in
the objective function to penalize the trace
\begin{eqnarray*}
&&\min \qquad  \bigl\langle\widetilde{\mathbf{X}}, \alpha \mathbf{I}_N-(1-
\lambda)\mathbf{A}+\lambda(\mathbf{J}_N -\mathbf {I}_N-
\mathbf{A}) \bigr\rangle
\\
&&\mbox{subject to}\qquad \widetilde{\mathbf{X}}\succeq\mathbf{0},
\\
&& \hspace*{20pt}\qquad\qquad 0\leq\widetilde{X}_{ij}\leq1 \qquad \mbox{for } 1\leq i, j\leq N,
\end{eqnarray*}
which is equivalent to
%
\begin{eqnarray}
\nonumber
&&\min  \qquad \langle\widetilde{
\mathbf{X}}, \mathbf {E} \rangle
\\
\label{eqcvx}
&&\mbox{subject to} \qquad  \widetilde{\mathbf{X}}\succeq\mathbf{0},
\\
\nonumber
&&\hspace*{20pt}\qquad\qquad 0\leq\widetilde{X}_{ij}\leq1 \qquad \mbox{for } 1\leq i, j\leq N,
\end{eqnarray}
where
%
\begin{equation}
\label{eqEdef}
\mathbf{E}:=\alpha\mathbf{I}_N-(1-\lambda)
\mathbf{A}+\lambda (\mathbf{J}_N -\mathbf{I}_N-
\mathbf{A}).
\end{equation}

\begin{remark}
At first glance, there are seemingly two tuning parameters: $\alpha$
and $\lambda$. In our theoretical result as shown later in Section~\ref{theory.sec}, the parameter $\alpha$ is required to be much greater
than the number of outlier nodes $m$. The introduction of $\alpha$
amounts to the trace penalization of $\widetilde{\mathbf{X}}$, which is
usually adopted in the literature of SDP relaxation in order to recover
a low-rank structure; see, for example, \citet{candes2012phaselift,li2012sparse}. In our problem, we intend to use (\ref{eqcvx}) to solve
for a low-rank matrix to reveal the clustering structure of the GSBM,
so this trace penalization is possibly a natural heuristic. However, in
our numerical simulations in Section~\ref{simulation.sec}, the
clustering effectiveness of the convex optimization method (\ref
{eqcvx}) is not significantly improved by choosing a positive $\alpha$.
Instead, (\ref{eqcvx}) works even by letting $\alpha$ be a small
constant or zero. On the contrary, there is a risk for choosing a large
$\alpha$, which may result in a positive definite $\mathbf{E}$. If
so, the
solution to (\ref{eqcvx}) must be $\mathbf{0}$, which is useless in
analyzing the networking data.

Therefore, we only need to tune the parameter $\lambda$ in practice,
and it has a clear statistical meaning as indicated in (\ref
{eqlambdapq}) in a special case of the ordinary SBM. In Section~\ref{theory.sec}, it is shown that if $\lambda$ lies in an interval
determined by $p^-$ and $q^+$ as defined in (\ref{eqdefp}), under
mild technical conditions, any solution $\widehat{\mathbf{X}}$ to
(\ref{eqcvx}) is capable of detecting the underlying group structure among
the inliers. A simple and heuristic data dependent choice of $\lambda$
is given in Section~\ref{simulation.sec}, where we also show
numerically that the performance of our method for clustering is robust
to the choice of $\lambda$.
\end{remark}

When $\widehat{\mathbf{X}}$ is obtained, in the pursuit of an explicit
clustering solution, a further step of $k$-means clustering is conducted
to the normalized column vectors of $\widehat{\mathbf{X}}$ with $k=r$,
provided the number of major clusters $r$ is assumed known.
Furthermore, in Section~\ref{theory.sec} it is shown that the
misclassification rate after the $k$-means clustering can be tightly controlled.

In summary, our proposed community detection procedure consists of the
following two steps:
\begin{longlist}[Step 2.]
\item[\textit{Step} 1.] Choose an\vspace*{1pt} appropriate tuning parameter $\lambda$, and
then solve (\ref{eqcvx}). The solution is denoted as $\widehat
{\mathbf{X}}$.

\item[\textit{Step} 2.]  Conduct $k$-means clustering algorithm to the normalized
column vectors of $\widehat{\mathbf{X}}$ with $k=r$, so that we can solve
for the assigning function $\hat{\phi}$ that maps from $\{1\leq i\leq
N\}$ to $\{1, \ldots, r\}$.
\end{longlist}

Finally, we introduce the augmented Lagrange multiplier algorithm
to
solve~(\ref{eqcvx}). Augmented Lagrange multiplier algorithms have
been employed in a variety of SDP optimizations in order to recover the
underlying low-rank matrix structure; see, for example, \citet{LLS2011,CLMW2011,JCSX2011,CSX2013} and a nice review paper on alternating
direction method of multipliers (ADMM) \citet{BPCPE2010}. Notice that
(\ref{eqcvx}) can be rewritten as
\begin{eqnarray*}
&&\min_{\mathbf{Y}, \mathbf{Z}}\qquad  \iota(\mathbf {Y}\succeq\mathbf{0})+\iota (
\mathbf{0} \leq\mathbf{Z}\leq\mathbf{J}_N)+ \langle\mathbf {Y},
\mathbf{E} \rangle,
\\
&&\mbox{subject to}\qquad \mathbf{Y}=\mathbf{Z},
\end{eqnarray*}
where the indicator function $\iota(a \in A)$ is defined as
\[
\iota(a \in A)= %
\cases{ 0, & $\quad a \in A$,
\vspace*{3pt}\cr
+\infty, & $\quad a \notin A$.}
\]
By this definition, we can easily conclude that $\iota(a \in A)$ is a
convex function if and only if $A$ is a convex set. Define the
augmented Lagrangian of this optimization problem as
\[
L_\rho(\mathbf{Y}, \mathbf{Z}; \bolds{\Lambda}):=\iota(\mathbf {Y}
\succeq\mathbf{0})+\iota(\mathbf{0} \leq\mathbf{Z}\leq\mathbf {J}_N)+
\langle\mathbf{Y}, \mathbf{E} \rangle+\frac{\rho}{2}\|\mathbf{Y}-\mathbf {Z}+
\bolds{\Lambda}\|_F^2.
\]

If both $\bolds{\Lambda}$ and $\mathbf{Z}$ are fixed, and we aim to minimize
$L_\rho(\mathbf{Y}, \mathbf{Z}; \bolds{\Lambda})$ with respect to~$\mathbf{Y}$,
it is equivalent to minimizing
\[
\iota(\mathbf{Y}\succeq\mathbf{0})+\frac{\rho}{2}\biggl\llVert \mathbf {Y}-
\mathbf{Z}+\bolds{\Lambda}+\frac{\mathbf{E}}{\rho}\biggr\rrVert _F^2.
\]
For any symmetric matrix $\mathbf{X}$ whose eigenvalue decomposition is
$\mathbf{V}\bolds{\Sigma}\mathbf{V}^{\intercal}$, define $\mathbf
{X}_+:=\mathbf{V}
\bolds{\Sigma}_+\mathbf{V}^{\intercal}$. Then the solution to the above
minimization has an explicit form
\[
\mathop{\operatorname{argmin}}_{\mathbf{Y}}L_\rho(\mathbf{Y}, \mathbf{Z}; \bolds {
\Lambda})= \biggl(\mathbf{Z}-\bolds{\Lambda}-\frac{\mathbf{E}}{\rho}
\biggr)_{+}.
\]

\begin{remark}
This step has dominating computational complexity in each iteration of
ADMM. In fact, an exact implementation of this subproblem of
optimization requires a full SVD of $\mathbf{Z}-\bolds{\Lambda
}-\frac{\mathbf{E}}{\rho}$, whose computational complexity is
$O(N^3)$. When $N$ is as
large as hundreds of thousands, the full SVD has scalability issue. An
open question is how to facilitate the implementation, or whether there
exists a surrogate that is computationally inexpensive. A possible
remedy is applying the low-rank iterative method, 
which means in each iteration of ADMM, the full SVD is replaced by a
partial SVD where only the leading eigenvalues and eigenvectors are
computed. Although this type of method may be stuck in local
minimizers, given the fact that SDP implementation can be viewed as a
preprocessing before $k$-means clustering, such a low-rank iterative
method might be helpful. We leave this large-scale computing problem as
a future research project.
\end{remark}

%

On the other hand, if both $\bolds{\Lambda}$ and $\mathbf{Y}$ are
fixed, to
minimize $L_\rho(\mathbf{Y}, \mathbf{Z}; \bolds{\Lambda})$ with
respect to
$\mathbf{Z}$ is equivalent to minimizing
\[
\iota(\mathbf{0}\leq\mathbf{Z}\leq\mathbf{J}_N)+\frac{\rho
}{2}
\llVert \mathbf{Z}-\mathbf{Y}- \bolds{\Lambda}\rrVert _F^2.
\]
Again, we have a closed-form solution
\[
\mathop{\operatorname{argmin}}_{\mathbf{Z}}L_\rho(\mathbf{Y}, \mathbf{Z}; \bolds {
\Lambda}):=\min \bigl(\max (\mathbf{Y}+\bolds{\Lambda}, \mathbf{0} ),
\mathbf{J}_N \bigr),
\]
which changes the negative entries of $\mathbf{Y}+\bolds{\Lambda}$ into
zeros and those greater than one into one.

As to the Lagrange multiplier, as the convention in the literature of
augmented Lagrange multiplier algorithms, $\bolds{\Lambda}$ is
updated to
$\bolds{\Lambda}+(\mathbf{Y}-\mathbf{Z})$.

\begin{algorithm}[b]
\caption{Robust community detection via alternating direction method}
\begin{algorithmic}
\label{algrcdadm1}
\STATE\textbf{Initialization:} $\mathbf{Z}_0=\mathbf{0}$, $\bolds
{\Lambda }_0=\mathbf{0}$, $\rho=1$ and $\mathit{iter}=100$.
\newline
\STATE\textbf{while $k<\mathit{iter}$}
\newline
\STATE1. $\mathbf{Y}_{k+1}:= (\mathbf{Z}_k-\bolds{\Lambda
}_k-\frac{\mathbf{E}}{\rho} )_{+}$;
\newline
\STATE2. $\mathbf{Z}_{k+1}:=\min (\max (\mathbf
{Y}_{k+1}+\bolds{\Lambda}_{k}, \mathbf{0}  ), \mathbf
{J}_N )$;
\newline
\STATE3. $\bolds{\Lambda}_{k+1}:=\bolds{\Lambda}_k+(\mathbf
{Y}_{k+1}-\mathbf{Z}_{k+1})$;
\newline
\STATE\textbf{end while}.
\newline
\STATE\textbf{Output the final} $\mathbf{Y}_{\mathit{iter}}$.
\end{algorithmic}
\end{algorithm}

The above augmented Lagrange multiplier method derives an iterative
algorithm for solving the convex optimization (\ref{eqcvx}), which is
summarized in Algorithm~\ref{algrcdadm1}. In numerical simulations,
we let $\mathbf{Z}_0=\mathbf{0}$ and $\bolds{\Lambda}_0=\mathbf
{0}$ for
initialization, and simply choose $\rho=1$ and run the algorithm for
$\mathit{iter}=100$ iterations. Numerical analyses of the algorithm applied to
simulated data and a real data set of political blogs are deferred to
Section~\ref{simulation.sec}, where its efficiency and effectiveness are clearly
demonstrated. Moreover, for the purpose of comparison, we also
implement ordinary spectral clustering methods on the synthetic data
sets. The numerical simulations clearly show that our method
outperforms spectral clustering methods in terms of robustness against outliers.


\section{Theoretical guarantees}
\label{theory.sec}
In this section, we will introduce our main theoretical results that
guarantee that the clustering procedure derived in the previous section
can detect the underlying communities under the GSBM. The following
theorem provides an explicit condition of the parameters $n$, $m$,
$p^-$, $q^+$ and $n_{\min}$, as well as the tuning parameters $(\alpha,
\lambda)$, under which the solution to (\ref{eqcvx}) is capable of
unveiling the underlying group structures among the inliers in presence
of a portion of outlier confounders.

\begin{thm}
\label{teomain}
Let $\mathbf{A}$ be the adjacency matrix of the semi-random graph under
the GSBM, as defined in \textup{(\ref{eqadjacency})}. Let $\widehat{\mathbf{X}}$
be a solution to the semidefinite program~\textup{(\ref{eqcvx})} and the
density gap $\delta$ be defined as in \textup{(\ref{eqdefdelta})}, and the
minimum within-group density $p^-$ and the maximum cross-group density
$q^+$ be defined as in \textup{(\ref{eqdefp})}. As defined in Section~\textup{\ref{intro.sec}}, the integer $n$ denotes the number of inlier nodes, $m$
denotes the number of outlier nodes and $n_{\min}$ denotes the minimum
community size among the inliers. Suppose that $p^-\geq C\frac{\log
n}{n_{\min}}$, $\alpha\geq3m$ and
%
\begin{equation}
\label{eqdeltalower0}
\delta>C \biggl(\sqrt{\frac{p^-\log n}{n_{\min}}}+\frac{\alpha
}{n_{\min
}}+
\frac{\sqrt{nq^+}}{n_{\min}}+\frac{m\sqrt{r}}{n_{\min}}+\frac
{nmp^-}{(\alpha- 2m)n_{\min}} \biggr)
\end{equation}
for some sufficiently large numerical constant $C$, and the tuning
parameter $\lambda$ satisfies
%
\begin{equation}
\label{eqlambdabound}
q^+ +\frac{\delta}{4}<\lambda<p^- -\frac{\delta}{4}.
\end{equation}
Then with probability at least $1-\frac{1}{n}-\frac{2r}{n^2}-\frac
{cr}{n_{\min}^4}$ for some numerical constant $c$, $\widehat{\mathbf{X}}$
must be of the form
%
\begin{equation}
\label{eqoutput}
\widehat{\mathbf{X}}=\mathbf{P}
\lleft[\matrix{
\mathbf{J}_{l_1} & & & \widehat{\mathbf {Z}}_1\vspace*{3pt}
\cr
&\ddots& & \vdots\vspace*{3pt}
\cr
& & \mathbf{J}_{l_r} & \widehat{
\mathbf {Z}}_r \vspace*{3pt}
\cr
\widehat {\mathbf{Z}}_1^{\intercal}
& \cdots& \widehat{\mathbf {Z}}_r^{\intercal} & \widehat{
\mathbf{W}} } \rright] %
\mathbf{P}^{\intercal},
\end{equation}
where $\mathbf{P}$ is defined as in \textup{(\ref{eqadjacency})}.
\end{thm}
%

Theorem~\ref{teomain}\vspace*{1.5pt} guarantees that any solution to (\ref{eqcvx})
$\widehat{\mathbf{X}}$ satisfies $\widehat{X}_{jk}=1$ for $\phi
(j)=\phi
(k)\leq r$, and $\widehat{X}_{jk}=0$ for $\phi(j)\neq\phi(k)$ and
$\phi
(j)\leq r, \phi(k)\leq r$. In other words, for each pair of inlier
nodes $j$ and $k$,\vspace*{1pt} whether they belong to the same group or not solely
depends on whether $\widehat{X}_{jk}$ equals $1$ or $0$. It is
noteworthy that condition~(\ref{eqlambdabound}) is similar to the
tuning parameter condition imposed in \citet{CSX2013}.

To interpret condition (\ref{eqdeltalower0}), it is helpful to
consider two examples. First, let us consider the very sparse case
where $p^- \simeq q^+ \simeq\delta\simeq O (\frac{\log
n}{n}
)$, $n_{\min} \simeq O(n)$ and hence $r \simeq O(1)$. This condition
implies that our procedure works even for a graph whose average degree
of inlier nodes is on the oder of $O(\log n)$. This is consistent with
the best-known result in the literature of community detection without
outliers by spectral clustering based on the adjacency matrices or
graph Laplacians [see \citet{LR2013}], although the $\log n$ barrier
could be resolved by more sophisticated nonbacktracking matrix
methods; see \citet{KMMNSZZ2013}. In this case, condition (\ref{eqdeltalower0}) becomes
\[
\delta> C \biggl(\frac{\log n}{n} + \frac{\alpha}{n} + \frac{m}{n} +
\frac{m}{\alpha- m}\frac{\log n}{n} \biggr).
\]
Then by letting $\alpha= \log N$, $m = \log n$ outliers are allowed.

In the second example, we assume $\delta\simeq p^- \simeq q^+ \simeq
O(1)$, and the number of clusters $r$ grows with $n$. As a specific
example, we let $r \simeq n^{{1}/{4}}$. Moreover, we assume
$n_{\min
} \simeq n^{{3}/{4}}$. Then condition (\ref{eqdeltalower0}) becomes
\[
1 \succsim\sqrt{\log n}n^{-{3}/{8}} + \alpha n^{-{3}/{4}} + \sqrt{\log
n}n^{- {1}/{4}} + m n^{-{5}/{8}} + n^{{1}/{4}} \frac{m}{\alpha- 2m}.
\]
Then by letting $\alpha= N^{{3}/{4}}$, $m=O(n^{{1}/{2} -
\varepsilon})$ outliers are allowed for any $\varepsilon> 0$.

%

A prominent feature of Theorem~\ref{teomain} is its consistency with
the state-of-the-art community detection under the ordinary SBM in the
literature. Assume there is no outlier node, that is, $m=0$, and we
simply let $\alpha=O(1)$ or just $\alpha=0$. Then condition (\ref{eqdeltalower0}) becomes
\[
\delta>C \biggl(\sqrt{\frac{p^-\log n}{n_{\min}}}+\frac{\sqrt{nq^+
\log
n}}{n_{\min}} \biggr).
\]
If the number of clusters is fixed, that is, $r=O(1)$, we also assume
the size of the smallest community $n_{\min}=O(n)$. As mentioned above,
this condition is guaranteed by letting the minimum within-group
density $p^-$ to be as low as $O(\frac{\log n}{n})$ and the density gap
$\delta=O(\frac{\log n}{n})$. In another example where $p^-=O(1)$,
$q^+=O(1)$ and $\delta=O(1)$, condition (\ref{eqdeltalower0}) is
equivalent to $n_{\min}\geq O(\sqrt{n \log n})$. By modifying Lemma~\ref
{teorandommatrix} as discussed in Section~\ref{proof.sec}, this
condition can be relaxed to $n_{\min}\geq O(\sqrt{n})$. This is
consistent with the state-of-the-art result in the community detection
literature by spectral clustering [see, e.g., \citet{CCT2012}], and
planted partition [see, e.g., \citet{GM2005,ST2007,OH2011,Ames2013,CSX2013}] where the within-group densities are usually assumed to be
the same, so do the cross-group densities. The $O(\sqrt{n})$ barrier of
the small cluster size is well known in the literature of planted
clique problems; see \citet{DM2013} and the references therein.

\begin{remark}
The proof of Theorem~\ref{teomain} is involved, and the details are
given in Section~\ref{proof.sec}. It is helpful to understand the
intuition behind the proof. The optimization (\ref{eqcvx}) consists of
two parts: a linear objective function and a constraint set which is
the intersection of a polytope and the semidefinite cone. In order to
show that the solution of (\ref{eqcvx}) has the form of (\ref
{eqoutput}), we find a point on the boundary of the constraint set such
that this point has the form of (\ref{eqoutput}). Moreover, we prove
that a level set of the linear objective function is tangent to the
tangent cone of the constraint set at the selected point. This shows
that the selected point is the solution of (\ref{eqcvx}). It is
noteworthy that the level set of the linear objective function is in
fact a hyperplane with co-dimension $1$, so the selected point is a
sharp vertex of the constraint set. For more details, see the remark
before the proof Lemma~\ref{thmsufficientcondition} in the
supplemental article \citet{CL2014}.
\end{remark}

Theorem~\ref{teomain} shows that the semidefinite programming (\ref
{eqcvx}) can discriminate the different groups among the inlier nodes.
However, the clustering result is not clear by only the observation of
$\widehat{\mathbf{X}}$, and it is not clear how the outliers could affect
the final clustering result. Given the extra knowledge of the number of
clusters, we propose to cluster the normalized column vectors of
$\widehat{\mathbf{X}}$ by $k$-means with parameter~$r$. To be specific,
without loss of generality, let us assume $\mathbf{P}=\mathbf{I}$,
and define
\[
\widehat{\mathbf{X}}=
\lleft[\matrix{ \mathbf{J}_{l_1}
& & & \widehat{\mathbf{Z}}_1 \vspace*{3pt}
\cr
&\ddots& & \vdots
\vspace*{3pt}
\cr
& & \mathbf{J}_{l_r} & \widehat{\mathbf
{Z}}_r\vspace*{3pt}
\cr
\widehat {\mathbf{Z}}_1^{\intercal}
& \cdots& \widehat{\mathbf {Z}}_r^{\intercal} & \widehat{
\mathbf{W}} } \rright]
=[\mathbf{x}_1, \ldots,
\mathbf{x}_N].
\]
Moreover, define $\mathbf{y}_i=\mathbf{x}_i/\|\mathbf{x}_i\|_2$.
Then all $\mathbf{y}_i$'s belong to the set of $N$-\break dimensional vectors
with two-norm $1$
and all coordinates being nonnegative. Notice that if \mbox{$\mathbf
{x}_i=\mathbf{0}$}, we then define $\mathbf{y}_i$ as an arbitrary
nonnegative vector
with norm $1$. Then, for any inlier indices $i, j\in I$ and $\phi
(i)\neq\phi(j)$, we have
\[
\|\mathbf{y}_i-\mathbf{y}_j\|_2^2=2-2
\mathbf{y}_i^{\intercal
}\mathbf{y}_j \geq 2-
\frac{2m}{n_{\min}},
\]
and for any $i, j\in I$ and $\phi(i)=\phi(j)=k$, we have
\[
\|\mathbf{y}_i-\mathbf{y}_j\|_2^2=2-2
\mathbf{y}_i^{\intercal
}\mathbf{y}_j \leq 2-
\frac{2 l_k}{l_k+m}=\frac{2m}{l_k+m}\leq\frac{2m}{l_k}.
\]
Moreover, for any $\mathbf{y}_i$ and $\mathbf{y}_j$, since both of
them are
nonnegative, we have
\[
\|\mathbf{y}_i-\mathbf{y}_j\|_2^2=2-2
\mathbf{y}_i^{\intercal
}\mathbf{y}_j \leq2.
\]
By definition, the solution to the $k$-means applied to $\{\mathbf{y}_1,
\ldots, \mathbf{y}_N\}$ is
%
\begin{equation}
\label{eqkmeans}
\mathop{\operatorname{argmin}}_{\mathcal{S}, \bolds{\mu}_1, \ldots, \bolds{\mu
}_r} \sum_{k=1}^r
\sum_{\mathbf{y}_j \in S_k}\|\mathbf{y}_j-\bolds{
\mu}_k\|^2,
\end{equation}
where $\mathcal{S}=\{S_1, \ldots, S_r\}$ is all $r$ nonoverlapping\vspace*{1pt}
partitions of $[N]$. It is obvious that $\bolds{\mu}_k = \frac
{1}{|S_k|}\sum_{\mathbf{y}_j \in S_k} \mathbf{y}_j$. We define
$D_i=\phi
^{-1}(k)$ for all $k=1,\ldots, r+1$, and choose $\tilde{\bolds{\mu}}_k$
as any vector $\mathbf{y}_i$ belonging to the $k$th community, that is,
$\phi(i)=k$. Then there holds
%
\begin{eqnarray}
\min_{\mathcal{S}, \bolds{\mu}_1, \ldots, \bolds{\mu}_r} \sum_{k=1}^r
\sum_{\mathbf{y}_j \in S_k}\|\mathbf{y}_j-\bolds{
\mu}_k\|^2 &\leq & \sum_{k=1}^{r-1}
\sum_{\mathbf{y}_j \in D_k}\|\mathbf{y}_j-\tilde{\bolds{
\mu}}_k\| ^2+\sum_{\mathbf{y}_j \in D_r \cup D_{r+1}}\|
\mathbf{y}_j-\tilde {\bolds{\mu}}_r\|^2
\nonumber
\\
&\leq & \sum_{k=1}^{r} \sum
_{\mathbf{y}_j \in D_k}\|\mathbf {y}_j-\tilde{\bolds{\mu
}}_k\|^2+\sum_{\mathbf{y}_j \in D_{r+1}}\|
\mathbf{y}_j-\tilde{\bolds{\mu}}_r\| ^2
\\
\nonumber
&\leq &  \Biggl(\sum_{k=1}^{r}
l_k\frac{2m}{l_k} \Biggr)+ 2m = 2mr+2m.
\end{eqnarray}

Suppose the solution to the $k$-means clustering is $\widehat{S}_1,
\ldots
, \widehat{S}_r$ and $\hat{\bolds{\mu}}_k = \frac{1}{|\widehat
{S}_k|}\sum_{\mathbf{y}_j \in\widehat{S}_k} \mathbf{y}_j$. For each $j \in
\widehat
{S}_k$, define $\hat{\phi}(j): =k$. Now we show that if $m < \frac
{n_{\min}}{2r+4}$, each $D_i$, $i=1, \ldots, r$ must account for more
than $50$ percent in some cluster $\widehat{S}_k$. Assume this is not
true. Then there is a $D_i$ being minority in each $\widehat{S}_k$, and
hence for each $\mathbf{y}_{a_j} \in D_i$, there exists a $\mathbf{y}_{b_j}
\notin D_i$, but $\hat{\phi}(\mathbf{y}_{a_j})=\hat{\phi}(\mathbf
{y}_{b_j})$.
Moreover, all these $2l_i$ indices are distinct. This implies
\[
\sum_{k=1}^r \sum
_{\mathbf{y}_j \in\widehat{S}_k}\|\mathbf {y}_j-\bolds{\hat{\mu
}}_k\|^2\geq\sum_{j=1}^{l_i}
\frac{1}{2}\| \mathbf{y}_{a_j}-\mathbf{y}_{b_j}\|
_2^2 \geq(l_i-m) \biggl(1-\frac{m}{n_{\min}}
\biggr).
\]
We then have $(n_{\min}-m) (1-\frac{m}{n_{\min}} )\leq
2m(r+1)$, which is contradictory to the assumption $m < \frac{n_{\min
}}{2r+4}$.

Since each $D_i$ is the majority of some estimated community $\widehat
{S}_k$, we can give the definition of misclassification rate among the
inliers: suppose there are $p$ pairs $(\mathbf{y}_{a_1}, \mathbf
{y}_{b_1}),\ldots, (\mathbf{y}_{a_p}, \mathbf{y}_{b_p})$ such that
all $2p$
indices are distinct, $1\leq\phi(\mathbf{y}_{a_j})< \phi(\mathbf
{y}_{b_j})\leq r$ for all $j=1, \ldots, p$ but $\hat{\phi}(\mathbf
{y}_{a_j})=\hat{\phi}(\mathbf{y}_{b_j})$. The misclassification rate among
the inliers is defined as $\max\frac{p}{n}$ for all possible $p$
satisfying the above property. Now we give an example showing why this
definition of misclassification rate is appropriate. Suppose $n$ balls
have $r$ colors as well as $m$ uncolored balls, and we assign them into
$r$ boxes. In the $i$th box, we assume there are $s_i$ balls having
color~$i$, while there are $t_i$ balls which are colored other than
$i$. Moreover, we also assume the assignment is acceptable in the sense
that $s_i>t_i$. In the $i$th box, there are at most $t_i$ distinct
pairs of colored balls such that in each\vspace*{1pt} pair the colors are different.
By our definition, the misclassification rate is $\frac{t_1+\cdots
+t_r}{n}$, which is the natural definition.

Back to our robust community detection problem, if we assume the
misclassification rate among the inliers is $\frac{p}{n}$, we have
\[
\sum_{k=1}^r \sum
_{\mathbf{y}_j \in\widehat{S}_k}\|\mathbf {y}_j-\bolds{\hat{\mu
}}_k\|^2\geq\frac{1}{2}\sum
_{j=1}^p \| \mathbf{y}_{a_j}-
\mathbf{y}_{b_j}\| _2^2 \geq p \biggl(1-
\frac{m}{n_{\min}} \biggr).
\]
Therefore, we have
\[
\frac{p}{n}\leq\frac{2mr+2m}{ (1-({m}/{n_{\min}})
)n}\leq \frac{(2r+3)m}{n}
\]
provided $m < \frac{n_{\min}}{2r+4}$. In summary, we have proven the
following theorem, which guarantees that the misclassification rate
among the inliers can be well controlled:
%
\begin{thm}
Suppose the assumptions in Theorem~\ref{teomain} hold as well as $m <
\frac{n_{\min}}{2r+4}$. Then, with high probability, the
misclassification rate among the inlier nodes $i \in I$ is less than or
equal to $\frac{(2r+3)m}{n}$.
\end{thm}

Rigorously speaking, $k$-means minimization is computationally NP-hard,
although in practice it is often easy and fast to implement with a
number of repetitions. However, as shown in \citeauthor{KSS2004}
[(\citeyear{KSS2004}),
Theorem~4.9], there is a $(1+\varepsilon)$ approximate $k$-means clustering
for (\ref{eqkmeans}) with computational time $O(2^{(r/\varepsilon
)^{O(1)}}N^2)$, which is polynomial time when $r$ is a constant.
Suppose $\{\check{S}_1, \ldots, \check{S}_r\}$ is a polynomial time
approximate $k$-means solution, such that
\begin{eqnarray*}
\sum_{k=1}^r \sum
_{\mathbf{y}_j \in\check{S}_k}\|\mathbf {y}_j-\bolds{\check{\mu
}}_k\|^2 &\leq & (1 + \varepsilon)\min_{\mathcal{S}, \bolds{\mu}_1, \ldots,
\bolds{\mu}_r}
\sum_{k=1}^r \sum
_{\mathbf{y}_j \in S_k}\|\mathbf {y}_j-\bolds{\mu }_k
\|^2
\\
&\leq & (1+\varepsilon) (2mr+2m).
\end{eqnarray*}
Then if within the inliers there are $p$ misclassified nodes by $\{
\check{S}_1, \ldots, \check{S}_r\}$, similarly to the previous argument, we
get $\frac{p}{n} \leq\frac{(1+\varepsilon)(2r+3)m}{n}$.

When $r$ grows with $N$, one can also cluster the rows of $\widehat
{\mathbf{X}}$ in (\ref{eqoutput}) based on the $\ell_1$ distance.
If two
inlier nodes belong to the same community, their corresponding rows in
$\widehat{\mathbf{X}}$ have $\ell_1$ distance less than $m$; on the other
hand, if two inlier nodes belong to different communities, their
corresponding rows have $\ell_1$ distance greater than $2n_{\min}$. If
the number of outliers is far less than the minimum size of the major
clusters, for example, $n_{\min} > O(m^2)$, a pairwise comparison
between the rows of $\widehat{\mathbf{X}}$ can detect the inlier
communities accurately even without the knowledge of~$r$. However, this
method does not work as effectively as $k$-means clustering in numerical
simulations. An interesting direction for future research is to figure
out whether there is a polynomial time $(1+\varepsilon)$ approximation
$k$-means clustering for~(\ref{eqkmeans}) when $r$ grows with $N$.


\section{Numerical results}
\label{simulation.sec}

In this section, synthetic data and a real-world network data are
employed to demonstrate the efficiency and effectiveness of our
community detection procedure: convex optimization (\ref{eqcvx})
followed by $k$-means. As discussed in Section~\ref{method.sec},
throughout all numerical simulations of the augmented Lagrange
multiplier method Algorithm~\ref{algrcdadm1}, we fix $\alpha=0$. All
simulations were carried out with MATLAB R2014b on a MacBook Pro with a
2.66 GHz Intel Core i7 Processor and 4GB 1067 MHz DDR3 memory. As
indicated in Algorithm~\ref{algrcdadm1}, for the initialization, let
$\bolds{\Lambda}_0=\mathbf{Z}_0=\mathbf{0}$. Also, we fix
$\mathit{iter}=100$ and $\rho
=1$. As to the $k$-means clustering to the normalized columns of
$\widehat
{\mathbf{X}}$, we use the ``kmeans'' function in MATLAB with $100$ replicates.

\subsection{Synthesized data simulations}
\label{synthetic.sec}

We consider again the synthetic data set used in Section~\ref{method.sec}. Figure~\ref{figbig2small2} illustrates the adjacency
matrix of a concrete realization of the original network. Suppose one
knows that there are $2$ major clusters, and a GSBM with $r=2$ clusters
is adopted.

We now explain in detail our implementation of Algorithm \ref{algrcdadm1}. First, we need to choose the tuning parameter $\lambda$
between the maximum cross-group density $q^+$ and the minimum
within-group density $p^{-}$. Ideal choices of $\lambda$ are formalized
by condition (\ref{eqlambdabound}) in Theorem~\ref{teomain}. In
practice, we propose a simple method to choose $\lambda$ as the mean
connectivity density in a subgraph determined by nodes with
``moderate''
degrees. If the adjacency matrix of the graph is denoted as $\mathbf{A}$,
and $\mathbf{1}_N$ represents the $N$-dimensional vector with all
coordinates equal to $1$, then $\mathbf{A}\mathbf{1}_N$ represents the
degrees of the $N$ nodes. The nodes with degrees above the $80$th
percentile or below the $20$th percentile are eliminated from the
graph, and $\lambda$ is chosen as the mean density of the \textit{subgraph} determined by the remaining nodes. The purpose of choosing
nodes with moderate degrees is to diminish the influence of the
outliers. Notice that the mean density of the subgraph may be very
different from the mean of $\mathbf{A}\mathbf{1}_N$, which is usually
significantly affected by the outliers.

\begin{figure}[t]

\includegraphics{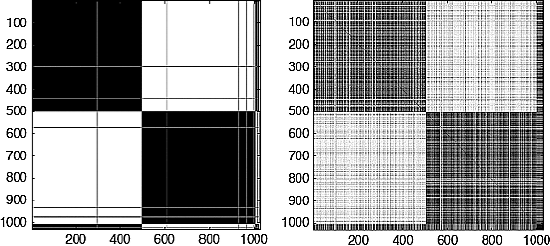}

\caption{On the right is the plot of the solution to convex
optimization (\protect\ref{eqcvx}). Based on it, the community detection
result followed by $k$-means algorithm is shown on the left.}
\label{figconvex}
\end{figure}

The convex method is implemented with $\lambda$ mentioned above. As an
illustration, in one realization of the synthetic data set, the
solution to convex optimization~(\ref{eqcvx}), and the community
detection result by further implementing $k$-means clustering with
$k=r=2$ are plotted in Figure~\ref{figconvex}.

We generated $10$ independent graphical data sets, and correspondingly
implemented $10$ trials of Algorithm \ref{algrcdadm1} followed by
$k$-means clustering, as well as spectral clustering on the graph
Laplacians and adjacency matrices. The average misclassification rate
among the $1000$ inlier nodes of our convex optimization method is
$0.0063$, which is much smaller than $1$ percent. The average time cost
for running Algorithm \ref{algrcdadm1} followed by $k$-means clustering
is $87.65$ seconds. In contrast, if we apply spectral clustering to the
graph Laplacians and adjacency matrices with $k=2$, respectively, the
average misclassification rates among the $1000$ inlier nodes are
$0.4792$ and $0.5000$, which are almost equivalent to random guessing.
If we treat the $30$ outliers as an additional group, and apply
spectral clustering to the graph Laplacians and adjacency matrices with
$k=3$, the misclassification rates among all $1030$ nodes are
correspondingly $0.3083$ and $0.4730$. Consequently, the
misclassification rates are high in terms of detecting the two major clusters.

\begin{figure}

\includegraphics{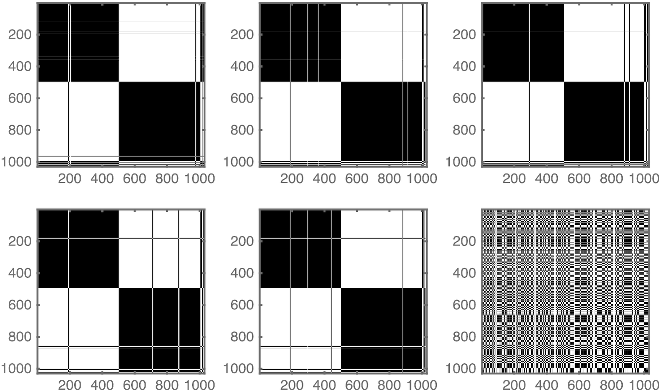}

\caption{The performance of convex community detection with different
values of $\lambda$.}
\label{figsensitivity}
\end{figure}

Now let us study the sensitivity of our algorithm to the choice of
$\lambda$. To be sure that $\lambda$ is between $q=0.11$ and $p=0.17$,
in Figure~\ref{figsensitivity} the community detection results are
plotted with $\lambda=0.11, 0.12, \ldots, 0.16$. It is obvious that for
our data set the clustering power is robust to $\lambda$, unless
$\lambda$ is too close to $p$. To our surprise, even when $\lambda=q$,
the two major clusters are well clustered. This is possibly due to the
facts that the graph is relatively sparse and the solution after $100$
iterations is still not exactly the solution to (\ref{eqcvx}).

\begin{figure}[b]

\includegraphics{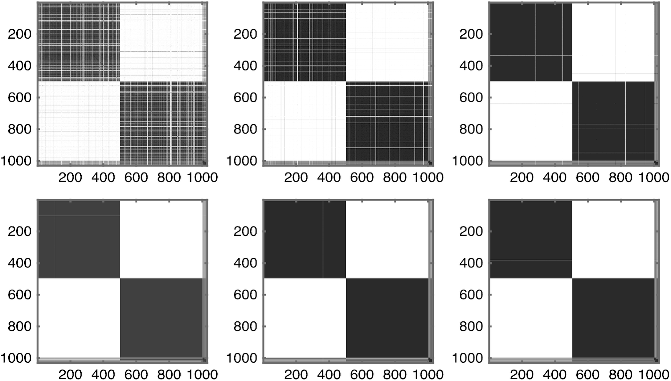}

\caption{The solutions of (\protect\ref{eqcvx}) with different
values of $p$.}
\label{figsolutions}
\end{figure}

On the right of Figure~\ref{figconvex}, we see that the solution to
(\ref{eqcvx}) is close to but not exactly equal to what Theorem~\ref
{teomain} predicts. A possible reason is that the density gap in our
synthetic data is not large enough. It is interesting that although the
solution does not have exactly the same form as in Theorem~\ref
{teomain}, the $k$-means in the second step can still successfully
cluster the two groups of nodes. We replace the within density $p=0.17$
with $0.19, 0.21, \ldots, 0.29$, and the solutions to (\ref{eqcvx})
are plotted in Figure~\ref{figsolutions}, respectively. The solutions
appear to be closer to the form in Theorem~\ref{teomain} as the
density gap increases.

\subsection{Real data application}
\label{application.sec}
%
\begin{figure}

\includegraphics{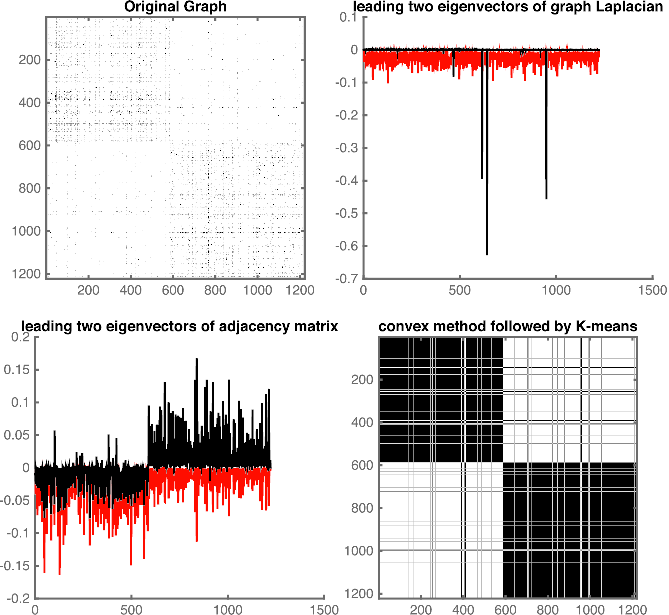}

\caption{Political blogs data of two clusters of conservatives and
liberals, along with the performance of convex optimization.}
\label{figplogosphere}
\end{figure}

In this section, our robust community detection procedure is tested by
implementing a modified version of convex optimization~(\ref{eqcvx})
on a political blogs network data set analyzed in \citet{AG2005}. This
network data set collected in 2005 is composed of political blogs and
their connections by hyperlinks, and it demonstrates the division and
interaction between the liberal and conservative blogs prior to the
2004 presidential election. By ignoring the directions of the
hyperlinks and selecting the largest connected component, there are
totally $1222$ nodes and $16{,}714$ edges, which implies that the
average degree is about $27$. As indicated in \citet{ZLZ2012}, the
distribution of the degrees is highly skewed to the right and has high
variability. Also, the political memberships of all blogs are clearly
studied and labeled manually in \citet{AG2005}, and are treated as the
truth for the purpose of evaluating the clustering efficacy of
different algorithms. The upper left panel of Figure~\ref{figplogosphere} plots the adjacency matrix of the observed political
blogs network.

Since the degrees in this real-world network data have high
variability, most community detection methods derived from the simple
SBM do not perform well. Instead, algorithms based on the so-called
degree-corrected SBM are proposed and proven to work well. For
instance, a polynomial time spectral method based on such a model is
introduced in \citet{CL2009}. Back to convex optimization (\ref
{eqcvx}), modification of the matrix $\mathbf{E}$ is needed to adapt to
the heterogeneity of the degrees. As mentioned earlier in Section~\ref{synthetic.sec} on the synthetic data simulation, $\lambda$ is chosen
data dependently as the mean of the degrees in a trimmed graph. When
the degrees have high variability, we propose to change the scalar
matrix $\lambda\mathbf{I}_N$ to the diagonal matrix $\mathbf{D}=\operatorname{Diag}(\mathbf{A}\mathbf{1}_N)/N$, the diagonal entries of which are the
degrees of all nodes divided by $N$. In brief, the modified convex
optimization is (\ref{eqcvx}) with
%
\begin{equation}
\label{eqdegreecorrectedE}
\mathbf{E}:=-(\mathbf{I}_N-\mathbf{D})^{{1}/{2}}
\mathbf {A}(\mathbf{I}_N-\mathbf{D})^{{1}/{2}}+
\mathbf{D}^{{1}/{2}}(\mathbf{J}_N -\mathbf{I}_N-
\mathbf{A})\mathbf{D}^{{1}/{2}}.
\end{equation}
In the second step of our proposed community detection procedure, we
choose $k=r=2$ in the $k$-means clustering. As a result, our community
detection procedure applied to the real-world network data set only
costs $137.16$ seconds to accurately cluster these $1222$ nodes with a
misclassification rate about $63/1222 \approx0.052$. The lower right
panel of Figure~\ref{figplogosphere} shows this clustering result by
plotting the adjacency matrix of the clustered graph, in which two
nodes are connected if and only if they are clustered in the same group.

The misclassification rate is comparable to the best-known results in
the literature. The SCORE method proposed in \citet{Jin2012} leads to a
misclassification rate of $58/1222$. Profile likelihood method under
degree-corrected SBM [\citet{KN2011}] and Newman--Girvan modularity
method [\citet{NG2004,ZLZ2012}] usually have misclassification rates
about $0.05$. However, as indicated in \citet{Jin2012}, the tabu
algorithm implemented to maximize these criteria is computationally
expensive and is numerically unstable due to bad initializations. It is
shown in \citet{Jin2012} that the average misclassification of the
modularity method is about $105/1222$ based on $100$ independent repetitions.

As to classical spectral clustering, the upper right and lower left
panels of Figure~\ref{figplogosphere} show that the two eigenvectors
of the graph Laplacian/adjacency matrix corresponding to the top two
eigenvalues are not capable in detecting and distinguishing the liberal
and conservative political blogs. Hence, ordinary spectral clustering
does not work when applied to this data set. A data-dependent penalized
spectral clustering applied to the graph Laplacian was proposed in
\citet
{JY2013}, but the misclassification rate is nearly $0.2$, which is much
worse than our result.

\section{Discussion}
\label{discussion.sec}

In this paper we introduce the GSBM for robust community detection in
the presence of arbitrary outlier nodes, and propose a computationally
feasible method using convex optimization. Strong theoretical
guarantees are established under mild technical conditions. In
particular, when the number of clusters is fixed and the edge density
within the inliers is $O({\log n\over n})$, $O(\log{n})$ outliers are
allowed; when the edge density within the inliers is on the order of
$O(1)$, and the number of clusters grows with $n$, for example,
$O(n^{{1}/{4}})$, our method is robust against $O(n^{{1}/{2} -
\varepsilon})$ adversarial outliers. Under the special case when there is
no outlier node, our theoretical result is also consistent with the
state-of-the-art results in the literature of computationally feasible
community detection under the SBM.
%

There are a number of possible extensions to the current results. The
proposed community detection procedure as well as the theoretical
guarantees depend on the assumption $\delta=p^- - q^+>0$. Although this
assumption is common in the literature of community detection, it is
actually a strong assumption which sometimes does not hold in
real-world network data applications. For example, suppose there are
$r=3$ clusters, and the connectivity matrix is
\[
\lleft[\matrix{ 0.4 & 0.2 & 0.05 \vspace*{3pt}
\cr
0.2 & 0.3 & 0.05
\vspace*{3pt}
\cr
0.05 & 0.05 & 0.1 } \rright].
\]
For each node, its associated within-group density is bigger than its
associated cross-group densities; however,
\[
\max_{1\leq i<j\leq r} B_{i j}>\min_{1\leq i\leq r}
B_{ii}.
\]
Therefore, in the current framework no choice of the tuning parameter
$\lambda$ is capable of the consistent community detection, which
implies the matrix $\mathbf{E}$ in the convex optimization step must be
modified. In fact, in our simulations, $\lambda$ is replaced by a
data-dependent diagonal matrix based on the degrees of all nodes in
order to adapt to high-degree variation. We are interested in
justifying this choice under the degree-corrected SBM proposed in \citet
{CL2009} and analyzed in \citet{KN2011,ZLZ2012,CCT2012,Jin2012,LR2013}.

In our numerical simulations, contrary to the established theoretical
guarantees, the choices of $\alpha$ are much smaller than the number of
outlier nodes $m$. In fact, the procedure works well with the choice
$\alpha=0$. An open question is whether this tuning parameter is
actually redundant. In addition, in the second step of our procedure,
the number of major inlier clusters $r$ is needed. Since the solution
of the convex optimization usually increases the connections within the
major groups and diminishes the connections across them, it is natural
and interesting to investigate whether $r$ can be inferred exactly from
the data. For reasons of space, we leave these as future\vspace*{-3pt} work.



\section{Proofs}
\label{proof.sec}

\subsection{Notation}
Throughout the proofs we will use the following notation: the $\ell\times
\ell$ identity matrix is denoted by $\mathbf{I}_{\ell}$. An $\ell_1
\times
\ell_2$ matrix whose entries all equal to $1$ is denoted as $\mathbf
{J}_{(\ell_1, \ell_2)}$. For square matrices, we write $\mathbf
{J}_{\ell
}:=\mathbf{J}_{(\ell, \ell)}$. An $\ell$-dimensional vector whose
coordinates all equal to $1$ is denoted as $\mathbf{1}_{\ell}$.

If all coordinates of a vector $\mathbf{v}$ are nonnegative, we write
$\mathbf{v}\geq\mathbf{0}$. When all coordinates of $\mathbf{v}$
are positive, we
write $\mathbf{v}>\mathbf{0}$.
We use $\mathbf{u}\geq\mathbf{v}$ to denote $\mathbf{u}-\mathbf
{v}\geq\mathbf{0}$,
and similarly, $\mathbf{u}>\mathbf{v}$ denotes $\mathbf{u}-\mathbf
{v} > \mathbf{0}$.
We also denote by $\|\mathbf{x}\|_\infty$ the maximum absolute values over
all coordinates of $\mathbf{x}$.

Similarly, if all entries of the matrix $\mathbf{M}$ are nonnegative, we
write $\mathbf{M}\geq\mathbf{0}$. When all entries of $\mathbf{M}$ are
positive, we write $\mathbf{M}>\mathbf{0}$. The inequality $\mathbf
{M}_1\geq\mathbf{M}_2$ denotes $\mathbf{M}_1-\mathbf{M}_2\geq
\mathbf{0}$, while $\mathbf{M}_1>\mathbf{M}_2$ denotes $\mathbf
{M}_1-\mathbf{M}_2 > \mathbf{0}$. Denote by $\|\mathbf{M}\|
_{\infty}$ the maximum absolute value over all entries of $\mathbf{M}$.
The norms $\|\cdot\|$ and $\|\cdot\|_F$ represent the operator and
Frobenius norms, respectively.

We use $\mathbf{M}\succ\mathbf{0}$ to denote that the symmetric
matrix $\mathbf{M}$ is positive definite and use $\mathbf{M}\succeq
\mathbf{0}$ to denote
that $\mathbf{M}$ is positive semidefinite. Similarly $\mathbf{M}_1
\succ\mathbf{M}_2$ and $\mathbf{M}_1 \succeq\mathbf{M}_2$
represent that $\mathbf{M}_1-\mathbf{M}_2$ is positive definite and
positive semidefinite, respectively.

For any vector $\mathbf{v} \in\mathbb{R}^n$, we denote by $\operatorname{Diag}(\mathbf{v})$ the $n\times n$ diagonal matrix whose diagonal
entries are
correspondingly the coordinates of $\mathbf{v}$.

Denote by $C, C_0, c$, etc. numerical constants, whose values could
change from line to line.

\subsection{Preliminaries}
Before proving Theorem~\ref{teomain}, we introduce several well-known
theorems in linear algebra and probability theory.

\begin{lemma}[{(Weyl  [\citet{HJ2013}, Theorem~4.3.1])}]\label{teoweyl}
Let $\mathbf{H}$ and $\mathbf{P}$ be two $n\times n$ Hermitian matrices.
Suppose that $\mathbf{H}+\mathbf{P}$, $\mathbf{H}$ and $\mathbf{P}$
have real eigenvalues
$\{\lambda_i(\mathbf{H}+\mathbf{P})\}_{i=1}^n$, $\{\lambda
_i(\mathbf{H})\}
_{i=1}^n$ and $\{\lambda_i(\mathbf{P})\}_{i=1}^n$, each arranged in
algebraically nonincreasing order. Then for $i=1, \ldots, n$ we have
\[
\lambda_i(\mathbf{H})+\lambda_n(\mathbf{P})\leq
\lambda_i(\mathbf {H}+\mathbf{P}) \leq\lambda_i(
\mathbf{H})+\lambda_1(\mathbf{P}).
\]
\end{lemma}

%
\begin{lemma}[{(Cauchy's interlacing theorem [\citet{HJ2013}, Theorem~4.3.28])}]
\label{teocauchy}
Let $\mathbf{H}$ be an $n\times n$ Hermitian matrix and $\mathbf{G}$
its $k
\times k$ principal submatrix. Suppose that $\mathbf{H}$ and $\mathbf{G}$
have real eigenvalues
$\{\lambda_i(\mathbf{H})\}_{i=1}^n$ and $\{\lambda_i(\mathbf{G})\}_{i=1}^k$,
each arranged in algebraically nonincreasing order. Then for $j=1,
\ldots, k$ we have
\[
\lambda_j(\mathbf{H})\geq\lambda_j(\mathbf{G})\geq
\lambda _{j+n-k}(\mathbf{H}).
\]
\end{lemma}

%
\begin{lemma}[{(Chernoff's inequality [\citet{Chernoff1981}])}]\label{teochernoff}
Let $X_1, \ldots, X_n$ be independent random variables with
\[
\mathbb{P}(X_i=1)=p_i, \qquad \mathbb{P}(X_i=0)=1-p_i.
\]
Then the sum $X=\sum_{i=1}^n X_i$ has expectation $\mathbb{E}(X)=\sum_{i=1}^n
p_i$, and we have
\begin{eqnarray*}
\mathbb{P}\bigl(X\leq\mathbb{E}(X)-\lambda\bigr)& \leq &  e^{-{\lambda
^2}/({2\mathbb{E}(X)})},
\\
\mathbb{P}\bigl(X\geq\mathbb{E}(X)+\lambda\bigr)& \leq &  e^{-{\lambda^2}/({2(\mathbb{E}(X)+\lambda/3)})}.
\end{eqnarray*}
\end{lemma}

Finally, we consider the following problem: suppose that $\mathbf
{A}=(a_{ij})_{1\leq i, j \leq n}$ is a random symmetric matrix, whose
diagonal entries are all zeros, while $a_{ij}, 1\leq i<j\leq n$ are
independent zero-mean Bernoulli random variables obeying $|a_{ij}|\leq
1$ and $\operatorname{Var}(a_{ij})\leq\sigma^2$. Can we prove that
with high
probability, $\|\mathbf{A}\|\leq C(\sigma\sqrt{n \log n}+ \log n)$ for
some numerical constant $C$? In the sequel, this upper bound is derived
by applying the following matrix Bernstein inequality, which is an
improvement of \citet{AW2002}:

%
\begin{lemma}[{[\citet{Tropp2011}, Theorem~6.1]}]
\label{teotropp}
Consider a finite sequence $\{\mathbf{X}_k\}$ of independent, random,
self-adjoint matrices with dimension $d$. Assume that
\[
\mathbb{E}\mathbf{X}_k = \mathbf{0} \quad \mbox{and} \quad \|
\mathbf{X}_k\|\leq R.
\]
If the norm of the total variance satisfies
\[
\biggl\llVert \sum_{k}\mathbb{E}\bigl(
\mathbf{X}_k^2\bigr)\biggr\rrVert \leq M^2,
\]
then the following inequality holds for all $t\geq0$:
\[
\mathbb{P} \biggl\{\biggl\llVert \sum_{k}
\mathbf{X}_k\biggr\rrVert \geq t \biggr\}\leq2d \exp \biggl(
\frac{-t^2/2}{M^2+Rt/3} \biggr).
\]
\end{lemma}

\begin{corollary}
\label{teowignerspec}
Let $\mathbf{A}=(a_{ij})_{1\leq i, j\leq n}$ be a symmetric random matrix
whose diagonal entries are all zeros. Moreover, suppose $a_{ij}$,
$1\leq i<j\leq n$ are independent zero-mean random variables satisfying
$|a_{ij}|\leq1$ and $\operatorname{Var}(a_{ij})\leq\sigma^2$. Then, with
probability at least $1-\frac{c}{n^4}$, we have
\[
\|\mathbf{A}\|\leq C_0 (\sigma\sqrt{n \log n}+ \log n )
\]
for some numerical constants $c$ and $C_0$.
\end{corollary}

\begin{pf}
For each pair $(i, j)\dvtx 1\leq i<j\leq n$, let $\mathbf{X}_{ij}$ be the
matrix whose $(i, j)$ and $(j, i)$ entries are both $a_{ij}$, whereas
other entires are zeros. Then we have
\[
\mathbf{A}=\sum_{1\leq i<j \leq n} \mathbf{X}_{ij}.
\]
Moreover, we can easily have $\mathbb{E}\mathbf{X}_{ij}=\mathbf{0}$,
$\|\mathbf{X}_{ij}\|\leq1$ and
\[
\mathbf{0}\preceq\sum_{1\leq i<j\leq n} \mathbb{E}\mathbf
{X}_{ij}^2 \preceq (n-1)\sigma^2
\mathbf{I}_n.
\]
They by applying Lemma~\ref{teotropp}, the proof is complete.
\end{pf}

\subsection{Supporting lemmas}
Notice that optimization (\ref{eqcvx}) is determined by the adjacency
matrix $\mathbf{A}$. Here we derive some properties of $\mathbf{A}$
and leave
the detailed proofs in the supplemental article \citet{CL2014}. More
precisely, we give some properties of the random matrix $\mathbf{K}$,
which is a principal submatrix of $\mathbf{A}$; see (\ref{eqadjacency}).

%
\begin{lemma}
\label{thmKaverage}
Recall that $p^-=\min_{1\leq i\leq r} B_{ii}$, $q^+=\max_{1\leq
i<j\leq
r} B_{i j}$ and $\delta=p^{-} - q^+$. If
%
\begin{equation}
\label{eqdeltalower}
\delta>C \biggl(\sqrt{\frac{q^+ \log n}{n_{\min}}}+\frac{\log
n}{n_{\min
}}
\biggr),
\end{equation}
for some sufficiently large numerical constant $C$, then with
probability at least $1-\frac{2}{n}-\frac{2r}{n^2}$, for all $i=1,
\ldots, r$ and $1\leq j<k\leq r$, we have
%
\begin{eqnarray}
\label{eqdiagrow}
\mathbf{K}_{ii}\mathbf{1}_{l_i} &\geq&
\bigl((l_i-1)B_{ii}-2\sqrt {(l_i-1)B_{ii}
\log n} \bigr)\mathbf{1}_{l_i},
\\
\label{eqoffrow}
\mathbf{K}_{jk}\mathbf{1}_{l_k} &\leq&
\biggl(B_{jk}+\frac{\delta
}{16} \biggr)l_k
\mathbf{1}_{l_j},
\\
\label{eqoffcolumn}
\mathbf{K}_{jk}^{\intercal}\mathbf{1}_{l_j}
&\leq& \biggl(B_{jk}+\frac{\delta
}{16} \biggr)l_j
\mathbf{1}_{l_k},
\\
\label{eqofftotal}
\mathbf{1}_{l_j}^{\intercal}\mathbf{K}_{jk}
\mathbf{1}_{l_k} &\geq& \biggl(B_{jk}-\frac{\delta}{16}
\biggr)l_kl_j.
\end{eqnarray}
\end{lemma}

%
\begin{lemma}
\label{teorandommatrix}
Suppose $p^-\geq C (\frac{\log n}{n_{\min}} )$. With
probability at least $1-c\frac{r}{n_{\min}^4}$, we have
%
\begin{equation}
\label{eqspectralnormdiagonal}
\bigl\llVert B_{ii} (\mathbf{J}_{l_i}-
\mathbf{I}_{l_i} )-\mathbf{K}_{ii}\bigr\rrVert \leq
C_0\sqrt{l_i B_{ii} \log l_i}, \qquad 1
\leq i\leq r
\end{equation}
and
%
\begin{equation}
\label{eqspectralnorm}
\|\mathbf{U}\|\leq C_0\bigl(\sqrt{nq^{+}
\log n}+ \log n \bigr),
\end{equation}
where $\mathbf{U}$ is an $n\times n$ symmetric matrix defined as
\[
\mathbf{U}:=
\lleft[ \matrix{ \mathbf{0} & \ldots&
B_{1r}\mathbf{J}_{(l_1, l_r)} - \mathbf{K}_{1r}
\vspace*{3pt}
\cr
\vdots& \ddots& \vdots \vspace*{3pt}
\cr
B_{1r}
\mathbf{J}_{(l_1, l_r)} - \mathbf{K}_{1r}^{\intercal} & \ldots &
\mathbf{0} } \rright]
\]
whose diagonal blocks are all $\mathbf{0}$'s. Here $C$, $C_0$ and $c$ are
some numerical constants.
\end{lemma}

It is worth noting that by applying a very recent result
\citeauthor{Vu2014} [(\citeyear{Vu2014}), Lemma~8], which is an improvement of \citet{FK1981,Vu2007},
we can
prove $\|\mathbf{U}\|\leq C_0(\sqrt{nq^{+}}+ \sqrt{\log n})$. Condition
(\ref{eqdeltalower0}) in Theorem~\ref{teomain} can then be relaxed to
\[
\delta>C \biggl(\sqrt{\frac{p^-\log n}{n_{\min}}}+\frac{\alpha
}{n_{\min
}}+\frac{\sqrt{nq^+}}{n_{\min}}+
\frac{m\sqrt{r}}{n_{\min}}+\frac
{nmp^-}{(\alpha- 2m)n_{\min}} \biggr).
\]
The benefit is that when $m=O(1)$, $p^-=O(1)$, $q^+=O(1)$ and $\delta
=O(1)$, $n_{\min}$ can be as small as $O(\sqrt{N})$ by letting
$\alpha
=\sqrt{N}$. In particular, if there is no outlier node, that is, the
ordinary SBM, this is consistent with the state-of-the-art result in
the literature of computationally feasible community detection.



\subsection{Proof of Theorem~\texorpdfstring{\protect\ref{teomain}}{3.1}}
In this section, we will rigorously prove Theorem~\ref{teomain}.
First, to simplify the calculations, we can assume the permutation
matrix $\mathbf{P}$ to be the identity matrix $\mathbf{I}_N$. This suggestion
is formalized by the following lemma:

\begin{lemma}
\label{teoequivalence}
If Theorem~\ref{teomain} is true for $\mathbf{P}=\mathbf{I}_N$, it
is also
true for any permutation matrix $\mathbf{P}$.
\end{lemma}

The proof is given in the supplemental article \citet{CL2014}. Lemma~\ref
{teoequivalence} guarantees that in order to prove Theorem~\ref
{teomain}, we can assume without loss of generality that $\mathbf
{P}=\mathbf{I}$, that is, $\mathbf{A}=
\bigl[{\fontsize{8.36}{10.36}{\selectfont \matrix{ \mathbf{K} & \mathbf{Z} \vspace*{0pt}\cr\mathbf
{Z}^{\intercal} & \mathbf{W}
} }}\bigr]$.

In the following, we will prove Theorem~\ref{teomain} based on the
following idea: In order to analyze a solution $\widehat{\mathbf{X}}$ to
(\ref{eqcvx}), we need to explore several inequalities that it
satisfies. The obvious ones are $\widehat{\mathbf{X}}\succeq\mathbf
{0}$ and
$\mathbf{0} \leq\widehat{\mathbf{X}} \leq\mathbf{J}_N$ as the feasibility
conditions in (\ref{eqcvx}). However,\vspace*{1pt} the optimality condition of
$\widehat{\mathbf{X}}$ implies that for any feasible $\widetilde
{\mathbf{X}}$, we have $\langle\widehat{\mathbf{X}}, \mathbf
{E}\rangle\leq\langle
\widetilde{\mathbf{X}}, \mathbf{E}\rangle$. To sufficiently utilize this
condition, we need to construct a feasible matrix $\mathbf{X}$, such that
$\langle\widehat{\mathbf{X}}, \mathbf{E}\rangle\leq\langle
\mathbf{X}, \mathbf{E}\rangle$ is a tight constraint. In Section~\ref{seccandidatesolution} we will show how to construct this $\mathbf{X}$.

After establishing these inequalities for any\vspace*{1.5pt} solution $\widehat
{\mathbf{X}}$, we give in Section~\ref{secsufficientcondition} a sufficient
condition which guarantees that $\widehat{\mathbf{X}}$ has the form
(\ref{eqoutput}) (with $\mathbf{P}=\mathbf{I}$), and then in
Section~\ref{secdualconstruction} we prove that with high probability this
sufficient condition is true by using the supporting lemmas proven
previously. Consequently, these three steps imply Theorem~\ref{teomain}.

\subsubsection{Solution candidate}

\label{seccandidatesolution}
In this section, we will construct a candidate solution $\mathbf{X}$
feasible to (\ref{eqcvx}). Denote
\begin{eqnarray*}
\mathbf{E}&=&\alpha\mathbf{I}_N+\lambda (\mathbf{J}_N-
\mathbf {I}_N )-\mathbf{A}
\\
&:=&
\lleft[\matrix{ (\alpha-\lambda)\mathbf{I}_{l_1}+
\lambda\mathbf{J}_{l_1}-\mathbf {K}_{11}& \ldots &\lambda
\mathbf{J}_{(l_1, l_r)}-\mathbf{K}_{1r}& \widetilde{
\mathbf{Z}}_1 \vspace*{3pt}
\cr
\vdots& \ddots& \vdots&\vdots
\vspace*{3pt}
\cr
\lambda\mathbf{J}_{(l_r, l_1)}-\mathbf{K}_{1r}^{\intercal}
& \ldots &(\alpha -\lambda)\mathbf{I}_{l_r}+\lambda\mathbf{J}_{l_r}-
\mathbf{K}_{rr}& \widetilde {\mathbf{Z}}_r \vspace*{3pt}
\cr
\widetilde{\mathbf{Z}}_1^{\intercal}& \ldots& \widetilde{\mathbf
{Z}}_r^{\intercal}& \widetilde{\mathbf{W}} } \rright],
\end{eqnarray*}
which is equivalent to defining
%
\begin{eqnarray}
\label{eqtildeZ} \widetilde{\mathbf{Z}}_i &=& \lambda
\mathbf{J}_{(l_i, m)}-\mathbf {Z}_i,\qquad  i=1, \ldots, r,
\\
\label{eqtildeW}
\widetilde{\mathbf{W}} &=& (\alpha-\lambda) \mathbf{I}_{m}+
\lambda \mathbf{J}_m -\mathbf{W}.
\end{eqnarray}
The following lemma, the proof of which is given in the supplemental
article \citet{CL2014}, guarantees the existence of $r$ vectors
$\mathbf{x}_1, \ldots, \mathbf{x}_r \in\mathbb{R}^m$, which will
be employed to
construct a candidate solution:

\begin{lemma}
\label{teoprimalconstruction}
If $\alpha\geq2m$ and $0<\lambda< 1$, the solution to
%
\begin{eqnarray}
\nonumber
&&\min \qquad \sum
_{i=1}^r \bigl\langle\mathbf{x}_i,
\widetilde{\mathbf {Z}}_i^{\intercal}\mathbf{1}_{l_i}
\bigr\rangle+\frac{1}{2}\sum_{i=1}^r
\mathbf{x}_i^{\intercal} \widetilde{\mathbf{W}}
\mathbf{x}_{i}
\\
\label{eqcvx2}
&& \mbox{subject to} \qquad \mathbf{x}_i\geq\mathbf{0}\qquad \mbox{for } 1
\leq i\leq r,
\\
&&\hspace*{20pt}\qquad\qquad\sum_{i=1}^r \mathbf{x}_i^{\intercal}
\bigl(\mathbf{e}_j\mathbf {e}_j^{\intercal
}\bigr)
\mathbf{x}_i\leq1 \qquad \mbox{for } 1\leq j\leq m,\nonumber
\end{eqnarray}
exists uniquely. Moreover, denote the solutions by $\mathbf{x}_1,
\ldots,
\mathbf{x}_r \in\mathbb{R}^m$, which by definition satisfy $\|
\mathbf{x}_i\|
_\infty\leq1$. Then there are nonnegative vectors
$\bolds{\beta}_1, \ldots, \bolds{\beta}_r \in\mathbb{R}^m$ and
an $m\times
m$ nonnegative diagonal matrix
\[
\bolds{\Xi}=\operatorname{diag}(\xi_1, \ldots, \xi_m),
\]
such that
%
\begin{eqnarray}
\label{eqxibeta}
\widetilde{\mathbf{W}}\mathbf{x}_i+\widetilde{
\mathbf {Z}}_i^{\intercal}\mathbf{1}_{l_i} &=& \bolds{
\beta}_i-\bolds{\Xi }\mathbf{x}_i,
\\
\label{eqxi}
\xi_j \Biggl(1- \sum_{i=1}^r
\mathbf{x}_i^{\intercal}\bigl(\mathbf {e}_j
\mathbf{e}_j^{\intercal}\bigr)\mathbf{x}_i \Biggr) &=& 0,\qquad
j=1,\ldots, m
\end{eqnarray}
and
%
\begin{equation}
\label{eqbeta}
\langle\mathbf{x}_i, \bolds{\beta}_i
\rangle=0,\qquad i=1, \ldots, r.
\end{equation}
For all $1\leq j, k\leq r$, there holds
%
\begin{equation}
\label{eqxixupper}
\mathbf{x}_j^{\intercal}(\widetilde{\mathbf{W}}+
\bolds{\Xi })\mathbf{x}_k \leq m \sqrt{l_jl_k}.
\end{equation}
Furthermore, for all $i=1, \ldots, r$ and $j=1, \ldots, m$, we have
%
\begin{equation}
\label{eqxibetainequality}
\beta_{i_j}+\mathbf{e}_j^{\intercal}
\mathbf{Z}_i^{\intercal
}\mathbf{1}_{l_i} \leq(\alpha-
\lambda+ \xi_j)x_{i_j}+\lambda l_i + \lambda
\sum_{k=1}^m x_{i_k}.
\end{equation}
Finally, for all $i=1, \ldots, r$,
%
\begin{equation}
\label{eqbetaupper}
\mathbf{0} \leq\bolds{\beta}_i \leq(m+l_i
- 1)\mathbf{1}_m.
\end{equation}
\end{lemma}

Throughout the paper, we define
\[
\mathbf{V}:=[\mathbf{v}_1, \ldots, \mathbf{v}_r]:=
\lleft[\matrix{ \mathbf{1}_{l_1} & \mathbf{0} & \ldots&
\mathbf{0} \vspace*{3pt}
\cr
\mathbf{0} & \mathbf{1}_{l_2} & \ldots&
\mathbf{0} \vspace*{3pt}
\cr
\vdots& \vdots& \ddots& \vdots \vspace*{3pt}
\cr
\mathbf{0} & \mathbf{0}& \ldots& \mathbf{1}_{l_r} \vspace*{3pt}
\cr
\mathbf{x}_1& \mathbf{x}_2 & \ldots&
\mathbf{x}_r } \rright] %
\]
and
\[
\mathbf{X}=\mathbf{V}\mathbf{V}^{\intercal}= %
\lleft[\matrix{
\mathbf{J}_{l_1} &\ldots& \mathbf{0} & \mathbf {1}_{l_1}
\mathbf{x}_1^{\intercal}\vspace*{3pt}
\cr
\vdots&\ddots& \vdots&
\vdots\vspace*{3pt}
\cr
\mathbf{0} & \ldots& \mathbf{J}_{l_r} &
\mathbf{1}_{l_r}\mathbf{x}_r^{\intercal
} \vspace*{3pt}
\cr
\mathbf{x}_1\mathbf{1}_{l_1}^{\intercal} & \cdots&
\mathbf{x}_r\mathbf{1}_{l_r}^{\intercal} & \mathbf
{x}_1\mathbf{x}_1^{\intercal}+\cdots+
\mathbf{x}_r\mathbf {x}_r^{\intercal}} \rright].
\]
Since $\mathbf{x}_i$'s are feasible to optimization (\ref{eqcvx2}), we
can easily see that $\mathbf{X}$ is feasible to optimization (\ref
{eqcvx}). We aim to prove that under mild technical conditions,
$\mathbf{X}$ is actually a solution to optimization (\ref{eqcvx}).

\subsubsection{Sufficient condition for the optimality of $\mathbf{X}$}
\label{secsufficientcondition}
In this section,\vspace*{1pt} we propose a condition which guarantees that any
solution $\widehat{\mathbf{X}}$ to (\ref{eqcvx}) must be in the form
of~(\ref{eqoutput}) with $\mathbf{P}=\mathbf{I}_N$. This sufficient condition
is equivalent to constructing a matrix $\bolds{\Lambda}$ satisfying a
series of equalities and inequalities as indicated in the following
lemma. We call it a dual certificate. In Section~\ref{secdualconstruction}, we will show that with high probability, this
dual certificate can be constructed in an explicit way.

\begin{lemma}
\label{thmsufficientcondition}
Suppose\vspace*{1pt} $\bolds{\Xi}$ and $\bolds{\beta}_1, \ldots, \bolds
{\beta}_r$ are
defined as in Lemma~\ref{teoprimalconstruction}. If there exist
symmetric matrices $\bolds{\Lambda}\in\mathbb{R}^{N\times N}$,
$\bolds{\Psi }_{jj}\in\mathbb{R}^{l_j\times l_j}$ $(1\leq j \leq r)$
and\vspace*{1pt} matrices $\bolds{\Phi}_{jk}\in\mathbb{R}^{l_j\times l_k}$
$(1\leq j< k\leq r)$, such that
%
\begin{equation}\label{eqlambdapsiphi}
\qquad\bolds{\Lambda}=
{\fontsize{8}{10}{\selectfont \lleft[\matrix{\ds (\alpha-
\lambda)\mathbf{I}_{l_1}+\lambda\mathbf{J}_{l_1}-\mathbf
{K}_{11}+ \bolds{\Psi}_{11}& \ldots&\ds\lambda
\mathbf{J}_{(l_1,
l_r)}-\mathbf{K}_{1r}-\bolds{\Phi }_{1r}
& \ds\widetilde{\mathbf {Z}}_1-\frac{1}{l_1}\mathbf{1}_{l_1}
\bolds{\beta }_1^{\intercal} \vspace*{3pt}
\cr
\vdots& \ddots&
\vdots&\vdots \vspace*{3pt}
\cr
\ds\lambda\mathbf{J}_{l_r, l_1}-
\mathbf{K}_{1r}^{\intercal}-\bolds {\Phi }_{1r}^{\intercal}
& \ldots&\ds(\alpha-\lambda)\mathbf {I}_{l_r}+\lambda\mathbf{J}_{l_r}-
\mathbf{K}_{rr}+\bolds{\Psi }_{rr} &\ds \widetilde{
\mathbf{Z}}_r-\frac
{1}{l_r}\mathbf{1}_{l_r}\bolds{
\beta}_r^{\intercal} \vspace*{3pt}
\cr
\ds\widetilde{
\mathbf{Z}}_1^{\intercal}-\frac{1}{l_1}\bolds{\beta
}_1\mathbf{1}_{l_1}^{\intercal} & \ldots& \ds\widetilde{
\mathbf {Z}}_r^{\intercal}-\frac
{1}{l_r}\bolds{
\beta}_r\mathbf{1}_{l_r}^{\intercal} & \ds\widetilde {
\mathbf{W}}+\bolds{\Xi}} \rright]}}\hspace*{-15pt}
\end{equation}
satisfies $\bolds{\Psi}_{ii}>\mathbf{0}$, $\bolds{\Phi
}_{jk}>\mathbf{0}$, $\bolds{\Lambda}\mathbf{V}=\mathbf{0}$ and
$\bolds{\Lambda}\succeq\mathbf{0}$,
%
then any minimizer $\widehat{\mathbf{X}}$ to (\ref{eqcvx}) must be
of the form
\[
\widehat{\mathbf{X}}= %
\lleft[\matrix{ \mathbf{J}_{l_1}
&\ldots& \mathbf{0} & \mathbf{1}_{l_1}\mathbf {x}_1^{\intercal
}+
\mathbf{H}_1 \vspace*{3pt}
\cr
\vdots&\ddots& \vdots& \vdots
\vspace*{3pt}
\cr
\mathbf{0} & \ldots& \mathbf{J}_{l_r} &
\mathbf{1}_{l_r}\mathbf {x}_r^{\intercal} +
\mathbf{H}_r \vspace*{3pt}
\cr
\mathbf{x}_1
\mathbf{1}_{l_1}^{\intercal} +\mathbf{H}_1^{\intercal}
& \cdots& \mathbf{x}_r\mathbf{1}_{l_r}^{\intercal}+
\mathbf{H}_r^{\intercal} & \mathbf{x}_1
\mathbf{x}_1^{\intercal}+\cdots+\mathbf{x}_r\mathbf
{x}_r^{\intercal}+\mathbf{H}_0 } \rright],
\]
which is the same as (\ref{eqoutput}). Moreover, $\mathbf{X}$ is a
solution to (\ref{eqcvx}).
\end{lemma}

An intuition behind the theorem and the rigorous proof are given in the
supplemental article \citet{CL2014}. It is noteworthy that the condition
on $\bolds{\Lambda}$ is weaker if the number of clusters $r$ gets
smaller. The reason is that the equality condition is $\bolds{\Lambda
}\mathbf{V}=\mathbf{0}$. Obviously when $r$ gets smaller, $\mathbf
{V}$ has fewer
columns, and hence the equality constraint becomes milder. We emphasize
that the choices of $\bolds{\Psi}_{ii}$ and $\bolds{\Phi}_{ij}$ are
intended to fit the equality constraint of $\bolds{\Lambda}$, that is,
$\bolds{\Lambda}\mathbf{V}=\mathbf{0}$. To make sure $\bolds
{\Lambda} \succeq\mathbf{0}$, we need to first project $\bolds
{\Lambda}$ onto the orthogonal
compliment of~$\mathbf{V}$, and then show the projection is positive
definite. This is based on the spectral norm bound as indicated in
Lemma~\ref{teorandommatrix}, which provides a concentration
inequality for a random matrix.

\subsubsection{Construction of dual certificate}
\label{secdualconstruction}
It suffices to construct a matrix $\bolds{\Lambda}$ in the form of
(\ref{eqlambdapsiphi}) in Lemma~\ref{thmsufficientcondition}, which
satisfies $\bolds{\Lambda}\mathbf{V}=\mathbf{0}$, $\bolds{\Psi
}_{ii}>\mathbf{0}$,
$\bolds{\Phi}_{jk}>\mathbf{0}$ and $\bolds{\Lambda}\succeq
\mathbf{0}$. The
following lemma guarantees the existence of such $\bolds{\Lambda}$, and
its proof is given in the supplemental article \citet{CL2014}.

%
\begin{lemma}
\label{teoLambdainequalities}
Suppose $p^-\geq C (\frac{\log n}{n_{\min}} )$, $q^+ +
\frac
{\delta}{4}<\lambda<p^- -\frac{\delta}{4}$ and $\alpha\geq3m$.
Moreover, assume
%
\begin{equation}
\label{eqdeltalower2} \qquad\delta>C \biggl(\sqrt{\frac{p^-\log n}{n_{\min}}}+\frac{\alpha
}{n_{\min
}}+
\frac{\sqrt{nq^+ \log n}}{n_{\min}}+\frac{m\sqrt{r}}{n_{\min
}}+\frac
{nmp^-}{(\alpha- 2m)n_{\min}} \biggr)
\end{equation}
for some sufficiently large numerical constant $C$. Then, with
probability at least $1-\frac{1}{n}-\frac{2r}{n^2}-\frac{cr}{n_{\min
}^4}$, there exist matrices $\bolds{\Psi}_{ii}\!$'s and $\bolds{\Phi
}_{jk}\!$'s
satisfying $\bolds{\Psi}_{ii}>\mathbf{0}$, $\bolds{\Phi
}_{jk}>\mathbf{0}$ and the
matrix $\bolds{\Lambda}$ defined by $\bolds{\Psi}_{ii}\!$'s and
$\bolds{\Phi }_{jk}\!$'s obey $\bolds{\Lambda}\mathbf{V}=\mathbf
{0}$ and $\bolds{\Lambda}\succeq
\mathbf{0}$.
\end{lemma}

\begin{supplement}[id=suppA]
\stitle{Supplemental materials to ``Robust and computationally
feasible community detection in the presence of arbitrary outliers nodes''\\}
\slink[doi]{10.1214/14-AOS1290SUPP} 
\sdatatype{.pdf}
\sfilename{aos1290\_supp.pdf}
\sdescription{We give in the supplement proofs to Lemmas \ref
{thmKaverage}, \ref{teorandommatrix}, \ref{teoequivalence}, \ref
{teoprimalconstruction}, \ref{thmsufficientcondition} and \ref
{teoLambdainequalities}.}
\end{supplement}

\printaddresses
\end{document}